\documentclass[11pt,reqno]{amsart}
\usepackage[utf8]{inputenc}
\usepackage[T1]{fontenc}
\usepackage{lmodern}
\usepackage[english]{babel}
\usepackage{amsmath,amssymb,amsthm,amsfonts,mathrsfs,mathtools}
\usepackage{microtype}
\usepackage{booktabs}
\usepackage{enumitem}
\usepackage{xcolor}
\usepackage{array}
\usepackage[margin=3cm, heightrounded]{geometry}
\usepackage{comment}
\usepackage[colorlinks=true,urlcolor=blue, citecolor=red,linkcolor=blue,
linktocpage,pdfpagelabels, bookmarksnumbered,bookmarksopen]{hyperref}
\usepackage[hyperpageref]{backref}
\mathtoolsset{showonlyrefs=true}
\hypersetup{
  pdftitle={Rotationally symmetric hypersurfaces with vanishing fractional mean curvature},
  pdfsubject={Rotationally symmetric stationary hypersurfaces for fractional mean curvature of order s in (0,1)},
  pdfkeywords={fractional mean curvature, rotational symmetry, stationary nonlocal hypersurfaces, concentric spheres, spectral stress, multiscale packing}
}
\allowdisplaybreaks
\numberwithin{equation}{section}

\theoremstyle{plain}
\newtheorem{theorem}{Theorem}[section]
\newtheorem{lemma}[theorem]{Lemma}
\newtheorem{proposition}[theorem]{Proposition}

\theoremstyle{definition}
\newtheorem{definition}[theorem]{Definition}

\theoremstyle{remark}
\newtheorem{remark}[theorem]{Remark}

\newcommand{\R}{\mathbb R}
\newcommand{\N}{\mathbb N}
\newcommand{\one}{\mathbf 1}
\newcommand{\dd}{\,d}
\newcommand{\TV}{\operatorname{Var}}
\newcommand{\kappaN}{\kappa_{n,s}}
\newcommand{\GammaN}{\Gamma_{n,s}}
\newcommand{\gammaN}{\gamma_{n,s}}
\newcommand{\An}{\mathfrak A_{n,s}}

\DeclareMathOperator{\dist}{dist}
\title{Rotationally symmetric fractional minimal surfaces}
\author[Gambicchia]{Chiara Gambicchia}
\address[C. Gambicchia]{
	Scuola Normale Superiore
	\newline\indent
	piazza dei Cavalieri 7,
	56127 Pisa, Italy}
\email{chiara.gambicchia@sns.it}
\author[Novaga]{Matteo Novaga}
\address[M.\ Novaga]{Dipartimento di Matematica
	\newline\indent
	Universit\`a di Pisa
	\newline\indent
	largo Bruno
Pontecorvo 5, 56127 Pisa, Italy}
\email{matteo.novaga@unipi.it}
\author[Talluri]{Matteo Talluri}
\address[M.\ Talluri]{Dipartimento di Matematica
	\newline\indent
	Alma Mater Studiorum  Universit\`a di Bologna
	\newline\indent
	piazza di Porta San Donato 5, 40126 Bologna, Italy}
\email{matteo.talluri@unibo.it}
\date{\today}


\begin{document}

\begin{abstract}
In this paper, we provide a non-trivial example of an $s$-minimal surface that is invariant under rotations.
More precisely, we study alternating radial sets whose boundaries are finite or locally finite unions of concentric spheres and prove that, fixing the outer radius to $r=1$, there exists a radial configuration for which the fractional mean curvature $H_s$ vanishes at every free interface.  
We then show that rescaling radial maximizers by their innermost radius and passing to the limit yields a properly embedded, rotationally symmetric, locally finite hypersurface satisfying $H_s=0$ at every point of the boundary.
\end{abstract}

\maketitle

\section{Introduction}\label{sec:intro}
The notion of fractional mean curvature and nonlocal minimal surfaces have been introduced  by Caffarelli, Roquejoffre and Savin in \cite{CRS}. 
It has been proved that bounded smooth sets with constant nonlocal mean curvature satisfy an Alexandrov-type rigidity theorem
\cite{CFSW,CFMN}; 
when the curvature vanishes, the sliding argument by S\'aez and Valdinoci in \cite{SV19} shows that $s$-minimal surfaces cannot be compact and, more generally, forces a smooth minimal set that is bounded in one direction to be a half-space. 
This result is even more investigated in \cite{CV}, where Cozzi and Valdinoci show that $s$-stationary sets are either half-spaces or have at least linear growth at infinity.
More recently, in \cite{CT}, Cozzi and Thompson prove a half-space theorem stating that any open set contained in a half-space, with boundary given by a smooth nonlocal minimal surface, is a half-space itself.
Another question of interest is the stability of $s$-minimal cones, which are proved to be flat in $\R^2$ for $s$ close to $0$ \cite{Cas25}, and in $\R^3$ for $s$ close to $1$ \cite{CCS}.
A very recent result is given in \cite{DSTV}, where the authors obtain optimal sheet separation estimates for stable $s$-minimal surfaces in dimension up to $7$ and prove a compactness theorem for stable $s$-minimal cones in $\R^8$.
On the other hand, periodic Delaunay-type hypersurfaces provide
a complementary unbounded constant-curvature picture \cite{CFW}.
Outside of the Euclidean setting, we refer the reader to \cite{CFS,CDSV} for the study of $s$-minimal surfaces on manifolds.
Finally, we refer the interested reader to \cite{Ser24} for a survey on nonlocal minimal surfaces and a discussion about future directions and applications.

However, apart from half-spaces, very few examples of $s$-minimal surfaces are known.
In \cite{DPW}, the reader can find a first non-trivial example of nonlocal minimal surface in $\R^3$ for $s\sim1$, the so-called fractional catenoid. Moreover, D\'avila, del Pino and Wei prove the existence of the minimal Lawson cones in the fractional setting, that is, they show that for any $m,n\geq 1$ and $s\in(0,1)$ there exists a value $\alpha$ such that the cone $C_{\alpha}=\{(u,v)\in \R^m\times\R^n\,:\,|v|\leq \alpha |u|\}$ is an $s$-minimal surface in $\R^{m+n}\setminus \{0\}$. Furthermore, in \cite{CDDp} the authors showed that helicoids have vanishing fractional mean curvature. 
In the present paper, we provide another nontrivial example of a non local minimal surface, which is given by an infinite union of concentric spheres. Unlike the known examples, our minimal surfaces has only bounded connected components. 
This kind of geometry arises naturally in fractional mean-curvature flow, for instance in the self-similar constructions of \cite{CN20}.

\subsection{Notation and main result}
Fix throughout $n\ge2$ and $s\in(0,1)$.  For $r>0$ let $x_r=re_1\in\partial B_r$.  Given
\begin{equation}\label{eq:radii}
 0<r_1<\cdots<r_M=1,\qquad M=2N,
\end{equation}
define the alternating radial set
\begin{equation}\label{eq:EN}
 E_N:=\bigcup_{k=1}^{N}\bigl(B_{r_{2k}}\setminus B_{r_{2k-1}}\bigr).
\end{equation}
At a regular boundary point we use the convention
\begin{equation}\label{eq:Hdef}
 H_s(x,E)
 :=\operatorname*{p.v.}\int_{\R^n}
 \frac{\one_{\R^n\setminus E}(y)-\one_E(y)}{|x-y|^{n+s}}\,\dd y.
\end{equation}
We call \eqref{eq:EN} a \emph{finite radial equilibrium} if
\begin{equation}\label{eq:equilibrium}
 H_s(x_{r_i},E_N)=0,\qquad i=1,\ldots,M-1.
\end{equation}
While we prove the existence of finite radial equilibria in Section~\ref{sec:existence-finite-equilibria} by a radial variational argument,
the classification and uniqueness of such equilibria remain open. 
The quantitative statements below apply to every equilibrium, not only to the variationally selected one. 
We remark that our finite configurations do not contradict the bounded rigidity results, since the fixed outer sphere is not required to satisfy the Euler--Lagrange equation.
The limiting object, where every interface is free and has zero fractional
mean curvature, is necessarily unbounded.

Before stating our main result, let us introduce the logarithmic and Euclidean gaps
\begin{equation}\label{eq:gaps}
 p_i:=\log\frac{r_{i+1}}{r_i}>0,
 \qquad
 g_i:=r_{i+1}-r_i>0,
 \qquad i=1,\ldots,M-1.
\end{equation}
Throughout the paper, the constants $c_{n,s},C_{n,s}$ may change from line to line and depend only on $n$ and $s$ unless further parameters are displayed.  
We remark that no estimate in this paper is claimed to be uniform as $s\downarrow0$ or $s\uparrow1$.

The main result of this paper is the following.
\begin{theorem}\label{thm:main}
There is a strictly increasing sequence
\begin{equation}\label{eq:infinite-radii-main}
 1=\rho_1<\rho_2<\cdots,
 \qquad \rho_i\longrightarrow+\infty,
\end{equation}
such that the set
\begin{equation}\label{eq:infinite-set-main}
 E_\infty:=\bigcup_{k=1}^\infty
 \bigl(B_{\rho_{2k}}\setminus B_{\rho_{2k-1}}\bigr)
\end{equation}
has locally finite fractional perimeter and $\partial E_\infty$
is a properly embedded smooth disconnected rotationally symmetric hypersurface
with
\begin{equation}\label{eq:infinite-stationary-main}
 H_s(x,E_\infty)=0 \qquad\text{for every }x\in\partial E_\infty
 .
\end{equation}
\end{theorem}
The above result is a direct consequence of the following two theorems, which provide finite equilibria for every even $M$ and the uniform geometric estimates needed to pass to the limit as $M\to+\infty.$
\begin{theorem}
    \label{thm:esistenza}For every even $M\ge2$ at least one finite radial equilibrium exists.
\end{theorem}
\begin{theorem}\label{thm:main finite}
Fix $n\ge2$ and $0<s<1$.  There are constants $c_{n,s},C_{n,s},\sigma_{n,s}>0$ such that every finite
radial equilibrium satisfies
\begin{equation}\label{eq:central-main}
 r_2-r_1\ge c_{n,s} r_1,
 \qquad\text{and hence}\qquad p_1\ge\log(1+ c_{n,s}).
\end{equation}
Moreover, there are $c_0,C_0,\eta_0>0$ such that
\begin{equation}\label{eq:adjacent-main}
 \min\{p_i,p_{i+1}\}\le\eta_0
 \quad\Longrightarrow\quad
 c_0p_i\le p_{i+1}\le C_0p_i,
 \qquad 1\le i\le M-3.
\end{equation}
Furthermore
\begin{equation}\label{eq:packing-main}
 r_1\le C_{n,s}M^{-\sigma_{n,s}}.
\end{equation}
Consequently $r_1\to0$ along every sequence of finite radial equilibria with
$M\to+\infty$.
\end{theorem}

\subsection{Plan of the paper}
In Section \ref{sec:preliminary}, we recall some preliminary results, which will be used throughout the paper, and we perform some preliminary computations in order to find an expression for the fractional curvature for a set of the form \eqref{eq:EN}. 
Then, in Section \ref{sec:existence-finite-equilibria}, we show the existence of finite radial equilibria, thus proving Theorem \ref{thm:esistenza}, and we provide an upper bound for the logarithmic gaps $p_i$, which will be needed to pass to the limit in the end.
In Section \ref{sec:technical-results} we prove a series of properties of the finite radial equilibria, which will be then used to prove Theorem \ref{thm:main finite}. In particular the lower bound for the logarithmic gaps will be used in Section \ref{sec:lower-bound-limit} to finally pass to the limit, thus concluding the proof of Theorem \ref{thm:main}.
\subsection*{Declaration of generative AI and AI-assisted technologies in the manuscript preparation process} 
The idea behind the construction of the minimal surface as a limit of shells is due to the authors. However, since the approximating shells are obtained through a maximization process and they are not explicit, it was very demanding to obtain estimates that allowed us to pass to the limit. These estimates, which are contained in Section \ref{sec:technical-results}, were obtained through several interactions with the model \textsf{GPT-6 Sol}. Furthermore, the model \textsf{Claude Opus 5.5} was used to check the correctness of the proofs and for the final editing. In this last check the model suggested a proof of Proposition \ref{prop:critical-identity} that is a lot simpler than the one obtained with the help of \textsf{GPT-6 Sol}. The first proof was based on a quite complicated approximation argument while the second one (that is, the one presented in this article) is a lot simpler and based only on the explicit value of the non-local curvature of the unit ball and the properties of the Pochhammer symbols.
All the proofs and the arguments have been checked by the authors, who take full responsibility for the content of the present article.

\section{Setting and preliminary results}\label{sec:preliminary}
In this section we recall some basic definitions and preliminary results, needed in the rest of the paper.
We start by recalling the definition of $s$-perimeter of a set.
\begin{definition}
For measurable sets $A,B\subset\R^n$, we consider the nonlocal interaction given by
\begin{equation}\label{eq:L-interaction}
 L(A,B):=\int_A\int_B\frac{\dd x\,\dd y}{|x-y|^{n+s}}.
\end{equation}
Then, the fractional perimeter of a set $A\subset\R^n$ is defined as
\begin{equation*}
    P_s(A):=L(A,A^c).
\end{equation*}
Moreover, following \cite{CRS}, given a domain $\Omega$ we set the fractional perimeter of $A$ in $\Omega$ as 
\[
P_s(A,\Omega) := L(A\cap\Omega,\,A^c\cap\Omega) + L(A\cap\Omega,\,A^c\setminus\Omega) + L(A\setminus\Omega,\,A^c\cap\Omega),
\]and we will say that $A$ has locally finite fractional perimeter if $P_s(A,B_R)$ is finite for every positive $R.$
\end{definition}
Thus, $2 P_s(A)$ is the homogeneous $W^{s,1}$ seminorm of
$\one_A$.
\begin{remark}
    From the definition, one can easily prove the following identities
    \begin{align}
      P_s(E\cup A)- P_s(E)
     &= P_s(A)-2L(A,E),
     &&\text{for all }A\subset E^c,
     \label{eq:add-toggle}\\
      P_s(E\setminus A)- P_s(E)
     &= P_s(A)-2L(A,E^c),
     &&\text{for all }A\subset E.
     \label{eq:remove-toggle}
    \end{align}
\end{remark}
\begin{definition}
    We say that $E$ is an $s$-minimal set in a bounded Lipschitz domain $\Omega$ if
    \[
    P_s(E,\Omega)\leq P_s(F,\Omega)
    \]
    whenever $E\setminus \Omega=F\setminus \Omega$.
\end{definition}
We now recall the notion of nonlocal curvature, as already introduced above.
\begin{definition}
Given a set $E\subset\R^n$ and a regular boundary point $x$, the $s$-mean curvature of $E$ at $x$ is given by
\begin{equation*}
 H_s(x,E)
 :=\operatorname*{p.v.}\int_{\R^n}
 \frac{\one_{\R^n\setminus E}(y)-\one_E(y)}{|x-y|^{n+s}}\,\dd y.
\end{equation*}
\end{definition}
In particular, analogously to the local case, if $E$ is an $s$-minimal set with regular boundary then
\[
H_s(x,E)=0
\]
at all $x\in\partial E$. Hence, again in analogy with the local case, we will say that $\partial E$ is an $s$-minimal surface if $H_s(x,E)=0$ for every $x\in\partial E.$

In order to compute many of the integrals which will appear in the paper, we remind the coarea formula.
\begin{theorem}[{\cite[Theorem 2.93]{AFP}}]
    Let $f:\R^m\to\R^k$ be a Lipschitz function and let $E$ be a countably $\mathcal{H}^n$-rectifiable set in $\R^m$, with $n\geq k$. 
    Then, it holds
    \[
    \int_EC_k d^Ef_x \dd \mathcal{H}^n(x)=\int_{\R^k}\mathcal{H}^{n-k}(E\cap f^{-1}(t))\dd t,
    \]
    where $C_k L:=\sqrt{\det(L\circ L^*)}$ is the $k$-dimensional coarea factor, and by $d^Ef_x$ we denote the tangential differential.
\end{theorem}

Finally, we recall the definitions and properties of some special functions, starting with Euler's Gamma function.
We refer to \cite{AAR} for all the details.
\begin{definition}
    For a complex number $z$, with Re$(z)>0$, we define the Gamma function at $z$ as the improper integral
    \[
    \Gamma(z)=\int_0^{+\infty}t^{z-1}e^{-t}\dd t.
    \]
\end{definition}
We remind that for all positive integer numbers $n$, it holds $\Gamma(n)=(n-1)!$ and that, in general, it holds $\Gamma(z+1)=z\Gamma(z)$.
Then, we can define Euler's Beta function as follows.
\begin{definition}
    Given $x,y$ complex numbers, with Re$(x)>0$, Re$(y)>0$, we define the Beta function at $(x,y)$ as the integral
    \[
    B(x,y)=\int_0^{1}t^{x-1}(1-t)^{y-1}\dd t.
    \]
\end{definition}
It is possible to prove, see \cite[Theorem 1.1.4]{AAR}, that the following identity holds
\[
B(x,y)=\frac{\Gamma(x)\Gamma(y)}{\Gamma(x+y)}.
\]
We now introduce the so-called Pochhammer symbols (also known as shifted factorial), defined as follows.
\begin{definition}
    For $m\in\N$ and $z\in\mathbb C$ with Re$(z)>0$, we write
    \[
     (z)_0:=1,
     \qquad
     (z)_m:=z(z+1)\cdots (z+m-1)=\frac{\Gamma(z+m)}{\Gamma(z)}.
    \]
\end{definition}
Lastly, we recall the definition of the function ${}_2F_1$.
\begin{definition}
    We set
    \[
    {}_2F_1(a,b;c;x):=\sum_{n=0}^{\infty}\frac{(a)_n(b)_n}{(c)_n}\,\frac{x^n}{n!}
    \]
\end{definition}
By \cite[Theorem 2.1.1]{AAR}, we know that the series on the right-hand side converges absolutely for $|x|<1$.

\subsection{Radial kernels and equilibrium equations}\label{sec:kernels}
In this section, we perform some preliminary computations, which will allow us to restate the equation \eqref{eq:equilibrium} in a more manageable way.
For $0<\rho<r<R$, set
\begin{equation}\label{eq:IJ}
 I(r,\rho):=\int_{B_\rho}\frac{\dd y}{|x_r-y|^{n+s}}\,,
 \qquad
 J(r,R):=\int_{\R^n\setminus B_R}\frac{\dd y}{|x_r-y|^{n+s}}\,,
\end{equation}
where $x_r=re_1$ as above. 
By an immediate change of variable $y=rz$, one can show that
\begin{equation*}
 I(r,\rho)=r^{-s}\mathcal I(\rho/r),
 \qquad
 J(r,R)=r^{-s}\mathcal J(R/r),
\end{equation*}
where we denote $\mathcal I(\rho/r):=I(1, \rho/r)$ and $\mathcal J(R/r):=J(1,R/r)$.
Finally, let
\begin{equation}\label{eq:kappa}
 \kappaN:=H_s(e_1,B_1).
\end{equation}
The following result provides the explicit value of of $\kappaN$ in terms of the data $n,s$.

\begin{lemma}\label{k esplicito}
The principal value in \eqref{eq:kappa} exists and is equal to
\begin{equation*}
\kappaN
 =\frac{2^{1-s}\pi^{(n-1)/2}\Gamma\bigl((1-s)/2\bigr)}
 {s\,\Gamma\bigl((n-s)/2\bigr)}>0.
\end{equation*}
\end{lemma}

\begin{proof}
For a smooth bounded set $E$ and $x\in\partial E$, one can see that
\[
 H_s(x,E)=\frac2s\,\operatorname*{p.v.}\!\int_{\partial E}
 \frac{(y-x)\cdot\nu_E(y)}{|x-y|^{n+s}}\,
 \dd\mathcal H^{n-1}(y),
\]
by applying the divergence theorem and using the fact that the flux across a sphere of radius $R$ is $O(R^{-s})$ and vanishes as $R\to+\infty$.
Hence, for $E=B_1$, $x=e_1$, and $y=\omega\in\mathbb S^{n-1}$, one has
$(\omega-e_1)\cdot\omega=1-\omega_1$ and
$|e_1-\omega|^2=2(1-\omega_1)$.
The coarea formula therefore yields
\begin{align*}
 \kappaN
 &=\frac{2^{1-(n+s)/2}}s|\mathbb S^{n-2}|
 \int_{-1}^1(1-u)^{1-(n+s)/2}(1-u^2)^{(n-3)/2}\,\dd u\\
 &=\frac{2^{-s}}s|\mathbb S^{n-2}|
 B\!\left(\frac{1-s}{2},\frac{n-1}{2}\right),
\end{align*}
which is exactly what we wanted to prove.
\end{proof}

Now, we have everything we need to exploit the curvature of a union of annuli $E_N$ with respect to a point $x\in \partial E_N$, as we show in the following result.
\begin{proposition}
Set $E_N$ as in \eqref{eq:EN}, and $M=2N$ as above. 
Then, for all indices $1\le i\le M$, the $s$-curvature of $E_N$ at a point $x_{r_i}\in\partial B_{r_i}$ can be written as
\begin{equation}\label{eq:curvature-formula}
 H_s(x_{r_i},E_N)
 =(-1)^i\frac{\kappaN}{r_i^s}
 +2\sum_{j<i}(-1)^{j+1}I(r_i,r_j)
 +2\sum_{j>i}(-1)^jJ(r_i,r_j).
\end{equation}
\end{proposition}

\begin{proof}
First of all, fix an index $i$ and let us write $\one_{E_N}$ and $\one_{\R^n\setminus E_N}$ as
\[
\one_{E_N}=\sum_{j=1}^{N}\one_{B_{r_{2j}}}-\one_{B_{r_{2j-1}}}
\]
and
\[
\one_{\R^n\setminus E_N}=1-\sum_{j=1}^{N}\one_{B_{r_{2j}}}-\one_{B_{r_{2j-1}}}.
\]
Hence, their difference is
\[
\one_{\R^n\setminus E_N}(y)-\one_{E_N}(y)=
1+2\sum_{j=1}^{M}(-1)^{j+1}\one_{B_{r_{j}}}(y).
\]
By reordering the terms and isolating the $i$-th term, we then obtain
\[
\one_{\R^n\setminus E_N}(y)-\one_{E_N}(y)=(-1)^i(1-2\one_{B_{r_i}}(y))+2\sum_{j=1}^{i-1}(-1)^{j+1}\one_{B_{r_j}}(y)+2\sum_{j=i+1}^M(-1)^{j}\one_{\R ^n \setminus B_{r_j}}(y),
\]
and the curvature at a point $x_{r_i}$ can be computed as
\begin{align*}
    H_s(x_{r_i},E_N)=&\text{p.v.}\int_{\R^n}\frac{\one_{\R^n\setminus E_N}(y)-\one_{E_N}(y)}{|x_{r_i}-y|^{n+s}}\dd y\\
    =&(-1)^i\,\text{p.v.}\int_{\R^n}\frac{1-2\, \one_{B_{r_i}}(y)}{|x_{r_i}-y|^{n+s}}\dd y+2\sum_{j=1}^{i-1}(-1)^{j+1}\int_{B_{r_j}}\frac{\dd y}{|x_{r_i}-y|^{n+s}}\\
    &+2\sum_{j=i+1}^M(-1)^{j}\int_{\R ^n\setminus B_{r_j}}\frac{\dd y}{|x_{r_i}-y|^{n+s}}\\
    =&(-1)^i\,H_s(x_{r_i}, B_{r_i})+2\sum_{j=1}^{i-1}(-1)^{j+1}I(r_i,r_j)+2\sum_{j=i+1}^M(-1)^{j}J(r_i, r_j).
\end{align*}
The thesis then follows simply noticing that $H_s(x_{r_i},B_{r_i})=r_i^{-s}H_s(e_1,B_1)=r_i^{-s}\kappaN$.
\end{proof}

Now, we define the radial density as follows
\begin{equation*}
 \mathcal K(r,\rho):=\int_{\partial B_\rho}\frac{1}{|x_r-y|^{n+s}}
 \,\dd\mathcal H^{n-1}(y).
\end{equation*}
This density allows us to estimate the functionals $I$ and $J$ when computed at $(r, r\pm d)$, for $d$ sufficiently small, as we can see in the following result.

\begin{lemma}\label{lem:near-contact}
Set
\begin{equation*}
 \gammaN:=\int_{\R^{n-1}}\frac{1}{(1+|z|^2)^{(n+s)/2}}\,\dd z
 \qquad\qquad \GammaN:=\frac{\gammaN}{s}.
\end{equation*}
Then, there are $\delta_0\in(0,1/8)$ and $C_{n,s}>0$ such that, for
$0<d\le\delta_0r$,
\begin{equation*}
 I(r,r-d)=\GammaN d^{-s}+O_{n,s}(r^{-s}),
 \qquad
 J(r,r+d)=\GammaN d^{-s}+O_{n,s}(r^{-s}),
\end{equation*}
while in general we have
\begin{equation}\label{IJ+-}
 I(r,r-d)\le C_{n,s}d^{-s}\quad \text{for $0<d<r$}\quad \text{and}\quad J(r,r+d)\le C_{n,s}d^{-s}\quad \text{for $d>0$}.
\end{equation}
\end{lemma}
\begin{proof}
By a straightforward change of variable one can easily see that
\[
\mathcal K(r,\rho)=\frac{1}{r^{1+s}}\mathcal K(1,\rho/r).
\]
Moreover, for all $r$ we have,
\begin{equation}\label{derivata}
 \partial_\rho I(r,\rho)=\mathcal K(r,\rho)>0
 \quad(0<\rho<r),
 \qquad
 \partial_RJ(r,R)=-\mathcal K(r,R)<0
 \quad(R>r).
\end{equation}
Indeed, (we show the computation for $I$, since for $J$ it is completely analogous) by the coarea formula it holds
\[
I(r,\rho)=\int_0^{\rho}\left(\int_{\partial B_{\sigma}}\frac{1}{|x_r-y|^{n+s}}\dd \mathcal{H}^{n-1}(y)\right)\dd \sigma=\int_0^{\rho}\mathcal{K}(r,\sigma)\dd\sigma
\]
which gives the derivative of $I$.

Therefore, it is enough to show that the following holds:
there exist $\delta_0,C_{n,s}>0$ such that, for $0<t\le\delta_0r$,
\begin{equation}\label{eq:density-expansion}
 \left|\mathcal K(r,r\pm t)-\gammaN t^{-1-s}\right|
 \le C_{n,s}r^{-1}t^{-s},
\end{equation}
and for all $r,\rho>0$, with $r\ne\rho$,
\begin{equation}\label{eq:Kglobal}
 \mathcal K(r,\rho)\le C_{n,s}\frac{1}{|r-\rho|^{1+s}}.
\end{equation}
Indeed, by \eqref{derivata} we have
\[
 I(r,r-d)=\int_d^r\mathcal K(r,r-t)\,\dd t,
 \qquad
 J(r,r+d)=\int_d^{+\infty}\mathcal K(r,r+t)\,\dd t.
\]
Hence, integrating \eqref{eq:density-expansion} where $d<t<\delta_0r$ and \eqref{eq:Kglobal} where $\delta_0r<t<r$ gives
\[
\left|I(r,r-d)-\int_d^r\frac{\gammaN}{t^{1+s}}\dd t\right|\leq C_{n,s}\int_d^r\frac{1}{rt^s}\dd t.
\]
Now notice that 
\[
 \gammaN\int_d^{\delta_0r}t^{-1-s}\,\dd t =\frac{\gammaN}{s}d^{-s}+O_{n,s}(r^{-s})  =\GammaN d^{-s}+O_{n,s}(r^{-s}),
\]
while the error gives
\[
 O\bigg(r^{-1}\int_d^{\delta_0r}t^{-s}\,\dd t\bigg) =O_{n,s}(r^{-s}),
\]
being $s\in(0,1)$.
The identity for $J$ is analogous. Now we observe that, by integrating \eqref{eq:Kglobal} we obtain, for any admissible $d$, \[
 I(r,r-d)\leq C_{n,s}\int_{d}^r\frac{dt}{t^{1+s}}\leq C_{n,s}d^{-s}
\]and \[
J(r,r+d)\leq C_{n,s}\int_{d}^{+\infty}\frac{dt}{t^{1+s}}\leq C_{n,s}d^{-s}
\]that is \eqref{IJ+-} and conclude the lemma.

So, let us now prove \eqref{eq:density-expansion} and \eqref{eq:Kglobal}.
By scaling, it is enough to consider $r=1$ and $\rho=1+\tau$, with $0<\tau<\delta_0$; the case with $\rho=r-\tau$ will be analogous.  
Decompose the unit sphere into a cap $A_1$ near $e_1$ and its complement $A_2$, clearly we have \[
\mathcal K(1,1+\tau)=\int_{(1+\tau)A_1}\frac{1}{|e_1-y|^{n+s}}d\mathcal H^{n-1}(y)+\int_{(1+\tau)A_2}\frac{1}{|e_1-y|^{n+s}}d\mathcal H^{n-1}(y).
\]We deal with the two integrals separately. The estimate for $A_2$ is simpler, indeed by a change of variable we have \[
\int_{(1+\tau)A_2}\frac{1}{|e_1-y|^{n+s}}d\mathcal H^{n-1}(y)=(1+\tau)^{n-1}\int_{A_2}\frac{1}{|e_1-(1+\tau)\omega|^{n+s}}d\mathcal H^{n-1}(\omega).
\]On $A_2$ we have $$|e_1-(1+\tau)\omega|\geq |e_1-\omega|-\tau\geq c_1$$uniformly for any $\tau$ small enough. Hence we find \[
\int_{(1+\tau)A_2}\frac{1}{|e_1-y|^{n+s}}d\mathcal H^{n-1}(y)=O_{n,s}(1)=O_{n,s}(\tau^{-s})\quad\text{as $\tau\to0.$}
\]
Now we deal with the integral on $(1+\tau)A_1.$ To do this we observe again that 
 \[
\int_{(1+\tau)A_1}\frac{1}{|e_1-y|^{n+s}}d\mathcal H^{n-1}(y)=(1+\tau)^{n-1}\int_{A_1}\frac{1}{|e_1-(1+\tau)\omega|^{n+s}}d\mathcal H^{n-1}(\omega).
\]The second integral can be estimated by considering the parametrization \[
\omega(z)=\big(\sqrt{1-|z|^2},z\big)\quad\text{for $z\in\R^{n-1}$ and $|z|\leq c_2$}
\]whose Jacobian is $(1-|z|^2)^{-\frac12}$. We thus obtain \[
\int_{(1+\tau)A_1}\frac{1}{|e_1-y|^{n+s}}d\mathcal H^{n-1}(y)=(1+\tau)^{n-1}\int_{|z|\leq c_2}\frac{dz}{|e_1-(1+\tau)\omega(z)|^{n+s}\sqrt{1-|z|^2}}.
\]By simple algebraic manipulations we find \[
|e_1-(1+\tau)\omega(z)|^2=\tau^2+(1+\tau)|z|^2+(1+\tau)\bigg(1-\sqrt{1-|z|^2}\bigg)^2,
\]that implies \[
\big||e_1-(1+\tau)\omega(z)|^2-(\tau^2+|z|^2)\big|\leq C(\tau|z|^2+|z|^4)\leq C_{\tau_0,c_2}|z|^2
\]for a constant that depends on $c_2$ and is uniform if $\tau\leq \tau_0$ is sufficiently small. Up to choosing $c_2$ and $\tau_0$ small enough, we find \[
\frac{1}{2}(\tau^2+|z|^2)\leq |e_1-(1+\tau)\omega(z)|^2\leq\frac{3}{2}(\tau^2+|z|^2).
\]By the above estimate and the mean value theorem we find \[
 \left||e_1-(1+\tau)\omega(z)|^{-n-s}-(\tau^2+|z|^2)^{-\frac{n+s}{2}}\right|\leq C_{n,s}(\tau+|z|^2)(\tau^2+|z|^2)^{-\frac{n+s}{2}}.
\]Observing that $$(1+\tau)^{n-1}(1-|z|^2)^{-\frac12}=1+O(\tau+|z|^2)$$ if $\tau$ and $c_2$ are sufficiently small, we find \begin{equation}\begin{aligned}\label{quattro integrali}
\bigg|\int_{(1+\tau) A_1}\frac{d\mathcal H^{n-1}(y)}{|e_1-y|^{n+s}}-\int _{|z|\leq c_2}\frac{dz}{(\tau^2+|z|^2)^{\frac{n+s}{2}}}\bigg|\leq &C_{n,s}\tau\int_{|z|\leq c_2}\frac{dz}{(\tau^2+|z|^2)^{\frac{n+s}{2}}}\\&+C_{n,s}\int_{|z|\leq c_2}\frac{|z|^2}{(|z|^2+\tau^2)^\frac{n+s}{2}}dz .\end{aligned}
\end{equation}By a change of variable it is easy to see that the two integrals on the right-hand side of \eqref{quattro integrali} are $O(\tau^{-s}),$ while the second one on the left-hand side can be estimated as \[
\int _{|z|\leq c_2}\frac{dz}{(\tau^2+|z|^2)^{\frac{n+s}{2}}}=\tau^{-1-s}\int_{\R^{n-1}}\frac{dz}{(1+|z|^2)^{\frac{n+s}{2}}}+O_{n,s}(1)
\]thus concluding the proof.

For \eqref{eq:Kglobal}, the same argument applies when $|r-\rho|\le cr$.  If $\rho\ge(1+c)r$, then $|x_r-y|\ge c'\rho$ on $\partial B_\rho$, so $\mathcal K(r,\rho)\le C\rho^{-1-s}\le C|r-\rho|^{-1-s}$. If $\rho\le(1-c)r$, then 
\[\mathcal K(r,\rho)\le C\rho^{n-1}r^{-n-s} \le Cr^{-1-s}\le C|r-\rho|^{-1-s}.\qedhere\]
\end{proof}
\begin{remark}
    The explicit value of $\gammaN$ is easily computable as 
    \[
    \gammaN=\pi^{(n-1)/2}\frac{\Gamma((1+s)/2)}{\Gamma((n+s)/2)},
    \]
    by using the polar coordinates and \cite[formula (1.1.20)]{AAR}.
\end{remark}

Before the next section, let us set some notation in order to simplify the rest of the proof.
Define
\begin{align}
 A_i^+&:=\sum_{k=1}^{M-i}(-1)^{k+1}J(r_i,r_{i+k})\label{A+i},\\
 A_i^-&:=\sum_{k=1}^{i-1}(-1)^{k+1}I(r_i,r_{i-k}).\label{A-i}
\end{align}
Then, \eqref{eq:equilibrium}, combined with \eqref{eq:curvature-formula} reads as
\begin{equation}\label{eq:balanced}
 A_i^+-A_i^-=\frac{\kappaN}{2r_i^s},
 \qquad \text{for all }i=1,\ldots,M-1.
\end{equation}
Notice that, since $I(r_i,r_j)$ is increasing in $r_j<r_i$ and $J(r_i,r_j)$ is decreasing in $r_j>r_i$, we have 
\[
0\leq A_i^+\leq J(r_i,r_{i+1})\qquad 0\leq A_i^-\leq I(r_i,r_{i-1}),
\]
for all $i=1,\dots , M-1$.


\section{Existence of finite radial maximizers}\label{sec:existence-finite-equilibria}
This section is devoted to the proof of Theorem \ref{thm:esistenza}.
In particular, here we prove the existence of finite radial configurations such that the equilibrium equation is satisfied at the internal radii.

Every finite radial configuration used in the variational problem
is bounded and has a finite union of smooth spherical interfaces, so its
fractional perimeter is finite.

For $M=2N$, define the closed simplex
\begin{equation*}
 \overline\Delta_M:=\{(r_1,\ldots,r_M):0\le r_1\le\cdots\le r_M=1\}.
\end{equation*}
For $\mathbf r\in\overline\Delta_M$ retain the definition
\begin{equation*}
 E(\mathbf r):=\bigcup_{k=1}^N \bigl(B_{r_{2k}}\setminus B_{r_{2k-1}}\bigr).
\end{equation*}
We remark that coincident radii are allowed in this subsection, and they simply cancel interfaces in pairs.
Let us begin with an estimate for the $s$-perimeter of an annulus.
\begin{lemma}\label{lem:thin-annuli}The following facts hold true:\begin{enumerate}[label=(\roman*)]
\item There exist a constant $C_{n,s}$ such that for any non negative $a$ and $\delta$ it holds \begin{equation}\label{eq:thin-annulus-two-sided1}
      P_s(B_{a+\delta}\setminus B_a)\leq C_{n,s}(a+\delta)^{n-1}\delta^{1-s}.
\end{equation}
\item For every $0<a_*<a^*<+\infty$ there are $c,\delta_*>0$ such that for all $a_*\le a\le a^*$ and $0<\delta<\delta_*$,
\begin{equation}\label{eq:thin-annulus-two-sided}
 c\delta^{1-s}\le  P_s(B_{a+\delta}\setminus B_a).
\end{equation}\end{enumerate}
\end{lemma}

\begin{proof}
Set $A=B_{a+\delta}\setminus B_a$.  For $x\in A$ set $d(x)=\dist(x,A^c)$.  
Then, since $A^c\subset\R^n\setminus B_{d(x)}(x)$, we have the following
\[
 \int_{A^c}\frac{1}{|x-y|^{n+s}}\,\dd y \le \int_{\R^n\setminus B_{d(x)}(x)}\frac{1}{|x-y|^{n+s}}\,dy=C_n\int_{d(x)}^{+\infty}\frac{1}{\rho^{s+1}}\dd \rho= C_{n,s}d(x)^{-s}.
\]
Since, $d(x)=\min\{|x|-a,a+\delta-|x|\}$, the polar coordinates give us
\begin{align*}
  P_s(A) \le C_{n,s}\int_a^{a+\delta} \rho^{n-1}\min\{\rho-a,a+\delta-\rho\}^{-s}\,\dd \rho \le C_{n,s}(a+\delta)^{n-1}\delta^{1-s},
\end{align*}
where $C$ depends on $n$ and $s$. This proves the upper bound in \eqref{eq:thin-annulus-two-sided}.

For the lower bound, we use the fact that, uniformly for $a\in[a_*,a^*]$, for $\delta$ sufficiently small, the spheres $\partial B_{a+\delta}$ satisfy a fixed exterior cone condition.  
Hence, if $x\in B_{a+\delta}\setminus B_{a+3\delta/4}$ and $t=a+\delta-|x|$, integrating on a cone contained in $B_{a+\delta}^c$, again via polar coordinates, we obtain
\[
 \int_{B_{a+\delta}^c}\frac{1}{|x-y|^{n+s}}\,\dd y \ge c_{n,s}t^{-s}.
\]
The area of every sphere with radius in $[a+3\delta/4,a+\delta]$ is bounded below by a positive constant depending only on $n,a_*$.  
Therefore,
\[
  P_s(A) \ge c\int_0^{\delta/4}t^{-s}\,\dd t = c'\delta^{1-s}.
\]
\end{proof}

\begin{lemma}\label{lem:perimeter-continuity}
The map
\[
 \mathbf r\longmapsto P_s(E(\mathbf r))
\]
is continuous on $\overline\Delta_M$. 
More precisely,
\begin{equation*}
 \left| P_s(E(\mathbf r))- P_s(E(\mathbf q))\right| \le C_{n,s,M}\sum_{i=1}^{M-1}|r_i-q_i|^{1-s}.
\end{equation*}
\end{lemma}

\begin{proof}
Reminding that $[\one_A]_{s,1}=2 P_s(A)$, the reverse triangle inequality for the Gagliardo seminorm gives
\[
 2\left| P_s(E(\mathbf r))- P_s(E(\mathbf q))\right| \le [\one_{E(\mathbf r)}-\one_{E(\mathbf q)}]_{s,1}.
\]
Now notice that
\[
 \one_{E(\mathbf r)} =\sum_{i=1}^M(-1)^i\one_{B_{r_i}}
\]
holds almost everywhere, also when two or more radii coincide.  
Hence, by the triangle inequality we have,
\begin{align*}
 [\one_{E(\mathbf r)}-\one_{E(\mathbf q)}]_{s,1} 
 \le\sum_{i=1}^{M-1}[\one_{B_{r_i}}-\one_{B_{q_i}}]_{s,1}
 =2\sum_{i=1}^{M-1} P_s\bigl(B_{\max\{r_i,q_i\}} \setminus B_{\min\{r_i,q_i\}}\bigr).
\end{align*}
Applying \eqref{eq:thin-annulus-two-sided1} in Lemma~\ref{lem:thin-annuli} we then conclude.
\end{proof}

The next lemma shows that the radial first variation of the perimeter coincides with the $s$-curvature, up to a multiplicative constant. The proof easily follows by writing the first variation of the fractional perimeter (see for instance \cite[Theorem 6.1]{F2M3}) with respect to the radial field $X(x)=\varphi(|x|)\,x/|x|$, where $\varphi\in C^\infty_c((r_{i-1},r_{i+1}))$ satisfies $\varphi\equiv1$ in a neighborhood of $r_i$, for any fixed index $i\leq M-1$. 
\begin{lemma}\label{lem:radial-first-variation}
Let $0<r_1<\cdots<r_M$ and let $E=E(\mathbf r)$.  For any fixed index $1\le i\le M-1$, set $r_i(t)=r_i+t$.  
Then, the following holds
\begin{equation}\label{eq:radial-first-variation}
 \left.\frac{\dd}{\dd t}\right|_{t=0} P_s(E(r_1,\ldots,r_i(t),\ldots,r_M)) =(-1)^i|\mathbb S^{n-1}|r_i^{n-1}H_s(x_{r_i},E).
\end{equation}
\end{lemma}
We are now ready to state and prove the existence of finite radial equilibria.
\begin{theorem}\label{thm:finite-existence}
For every even $M\ge2$, the function
$\mathbf r\mapsto P_s(E(\mathbf r))$ has a maximizer
$\mathbf r^{(M)}$ on $\overline\Delta_M$, and every maximizer satisfies
\begin{equation}\label{eq:maximizer-interior}
 0<r_1^{(M)}<\cdots<r_M^{(M)}=1.
\end{equation}
Consequently $E(\mathbf r^{(M)})$ is a finite radial equilibrium.
\end{theorem}

\begin{proof}
The existence of a maximizer follows from Lemma~\ref{lem:perimeter-continuity} and the compactness of $\overline\Delta_M$.  
We only need to prove that any maximizer is in the interior of $\Delta_M$.

First suppose that the only missing effective interface is the central one:$r_1=0<r_2<\cdots<r_M=1$.  
Then, $B_{r_2}\subset E$ up to a null set.
For $0<\delta<r_2/2$ let us consider $B_\delta$.  
Formula \eqref{eq:remove-toggle} and the positive distance between $B_{\delta}$ and $E^c$ give
\[
  P_s(E\setminus B_\delta)- P_s(E) = P_s(B_1)\delta^{n-s}-2L(B_\delta,E^c) \ge c_{n,s}\delta^{n-s}-C_{n,s,r_2}\delta^n,
\]
since we can easily estimate
\[
L(B_{\delta}, E^c)\leq \int_{B_{\delta}}\int_{B_{r_2}^c}\frac{1}{|x-y|^{n+s}}\dd x\dd y\leq C_{n,s,r_2}\,|B_{\delta}|.
\]
In particular, for $\delta$ small enough we can conclude that $P_s(E\setminus B_\delta)- P_s(E)>0$, against the maximality of $P_s(E)$.

Consider now an arbitrary remaining boundary point of $\Delta_M$.  
We can group equal positive radii into maximal consecutive blocks.  
A block of length $\ell$ at a radius $a>0$ contributes one effective spherical interface when $\ell$ is odd and no interface when $\ell$ is even. 
In either case, all unused indices in the block come in canceling consecutive pairs.
A block of zero radii contributes no geometric interface.

We first treat a terminal block
\[
 r_j=\cdots=r_M=1
\]
of even length $\ell=M-j+1\ge2$. 
Since $M$ and $\ell$ are even, $j$ is odd.
Replacing $r_j=1$ by $r_j=1-\delta$ and leaving the other radii fixed changes
the indicator by
\[
 -\one_{B_{1-\delta}}+\one_{B_1} = \one_{B_1\setminus B_{1-\delta}}.
\]
Thus, a thin outer shell $A_\delta=B_1\setminus B_{1-\delta}$ is added to $E$.  
Because the terminal block has even multiplicity, the phase immediately inside $\partial B_1$ is $E^c$; the nearest effective interface of $E$, if it exists, is a fixed positive distance away.  
Hence, similarly as above, we can estimate
\[
 L(A_\delta,E)\le C\delta,
\]
and \eqref{eq:add-toggle} together with Lemma~\ref{lem:thin-annuli} gives
\[
  P_s(E\cup A_\delta)- P_s(E) \ge c\delta^{1-s}-C\delta
\]
which is positive for small $\delta$, thus contradicting maximality.
Analogously, we obtain the same contradiction if the terminal block at $1$ has odd multiplicity or if the coinciding radii are strictly smaller than $1$.
Hence, \eqref{eq:maximizer-interior} holds.
Since the maximizer is in the interior with respect to $r_1,\ldots,r_{M-1}$, all the derivatives in Lemma~\ref{lem:radial-first-variation} vanish.  
The positive area factor then gives \eqref{eq:equilibrium}.
\end{proof}

\begin{remark}
The above Theorem maximizes $P_s$ only over the finite-dimensional radial class. 
It supplies an Euler--Lagrange equilibrium, not a minimizer under arbitrary perturbations.
Indeed, the construction gives no stability inequality outside the radial class.
\end{remark}

\subsection{Uniform upper bound for the logarithmic gaps}
In the last part of this section, we prove that the gaps $p_i$ all satisfy a uniform upper bound, which does not depend on the number of shells $M$.
Let us set
\[
p_i^{(M)}=\log \bigg(\frac{r_{i+1}^{(M)}}{r_{i}^{(M)}}\bigg),
\]
where the $r_i^{(M)}$'s are given by Theorem \ref{thm:finite-existence}.
\begin{lemma}\label{lem:gap-upper-all}
There is $\overline p_{n,s}<+\infty$ such that
\begin{equation*}
 p_i^{(M)}\le\overline p_{n,s},
 \qquad 1\le i\le M-1,
\end{equation*}
for every finite radial equilibrium.
\end{lemma}

\begin{proof}
For the sake of readability, we fix $M$ and omit it in the notation for the rest of the proof.
From \eqref{eq:balanced}, the positivity and first-term bound for alternating
sums give
\[
 J(r_i,r_{i+1})\ge A_i^+ = A_i^-+\frac{\kappaN}{2r_i^s} \ge \frac{\kappaN}{2r_i^s}.
\]
By scaling,
\[
 J(r_i,r_{i+1})=r_i^{-s}\mathcal J(e^{p_i}).
\]
Thus, $\mathcal J(e^{p_i})\ge\kappaN/2$, independently of $M$.  
Notice that the function $q\mapsto\mathcal J(q)=J(1,q)$ is continuous and strictly decreasing on $(1,+\infty)$ and tends to zero as $q\to+\infty$.
Hence, we can choose $\overline p_{n,s}$ so that $\mathcal J(e^{\overline p_{n,s}})<\kappaN/2$.
Then $p_i<\overline p_{n,s}$, and we have the desired bound.
\end{proof}
In order to be able to take this construction to the limit as $M\to +\infty$, we need to prove that a uniform lower bound holds as well. 
However, the proof of such bound will is much harder to obtain and will require some technical results, that we gather in the next section.

\section{Properties of finite radial equilibria}\label{sec:technical-results}
In this section we prove several properties of any finite radial equlibria. These results will be used in the proof of Theorem \ref{thm:main finite}.
\subsection{An asymptotic formula for the fractional curvature}
In this subsection we derive a power expansion for the interaction kernels $I$ and $J$ and then we derive an asymptotic formula for the fractional mean curvature with respect to the logarithmic gaps. 

For every $n\geq 2$ natural number and $s\in(0,1)$ we define
\begin{equation}\label{ans}
 \mathscr A_{n,s}(t):=\int_{\mathbb S^{n-1}}
 \frac{1}{(1-2t\omega_1+t^2)^{(n+s)/2}}\,\dd\omega,
 \qquad |t|<1.
\end{equation}where $\omega_1=\langle \omega;e_1\rangle.$

\begin{lemma}
\label{lem:angular}
For $|t|<1$ the function $\mathscr A_{n,s}$ can be expressed in power series as
\begin{equation}\label{eq:angular-series}
 \mathscr A_{n,s}(t)=\sum_{m=0}^{+\infty} a_m^{(n,s)}t^{2m},
 \qquad
 a_m^{(n,s)}=|\mathbb S^{n-1}|
 \frac{(\frac{n+s}{2})_m(1+\frac s2)_m}{(\frac n2)_m m!}>0.
\end{equation}
Moreover, the following asymptotic expansion holds true as $m\to+\infty.$
\begin{equation}\label{eq:am-asymptotic}
 a_m^{(n,s)}=\An m^s+O_{n,s}(m^{s-1}),
 \qquad
 \An=|\mathbb S^{n-1}|
 \frac{\Gamma(n/2)}{\Gamma((n+s)/2)\Gamma(1+s/2)}.
\end{equation}
\end{lemma}
\begin{proof}
    Since the integrand depends only on $\omega_1$ we have, by a direct computation via the coarea formula
    \[
\mathscr A_{n,s}(t)=\int_{-1}^1\dd\omega_1\int_{\mathbb S^{n-2}} \frac{(1-\omega_1^2)^{\frac{n-3}{2}}}{(1-2t\omega_1+t^2)^{\frac{n+s}2}}\dd \theta=|\mathbb S^{n-2}|\int_{-1}^1 \frac{(1-\omega_1^2)^{\frac{n-3}{2}}}{(1-2t\omega_1+t^2)^{\frac{n+s}2}}\dd\omega_1.
    \]
    By the change of variable $\omega_1=1-2x$ we find 
    \[
\mathscr A_{n,s}(t)=2^{n-2}|\mathbb S^{n-2}|(1-t)^{-n-s}\int_0^1x^\frac{n-3}{2}(1-x)^\frac{n-3}{2}\bigg(1+\frac{4tx}{(1-t)^2}\bigg)^{-\frac{n+s}{2}}\dd x. 
    \]
    By \cite[Theorem 2.2.1]{AAR} we can compute this last integral to obtain 
    \[
\int_0^1x^\frac{n-3}{2}(1-x)^\frac{n-3}{2}\bigg(1+\frac{4tx}{(1-t)^2}\bigg)^{-\frac{n+s}{2}}\dd x=B\bigg(\frac{n-1}2,\frac{n-1}{2}\bigg)(1-t)^{n+s}{}_2F_1\bigg(\frac{n+s}{2},1+\frac{s}{2};\frac{ n}{2};t^2\bigg)
    \]
    where we remind that $B$ denotes Euler's beta function. Thanks to the definition of the function ${}_2F_1$, we obtain \eqref{eq:angular-series}. To prove \eqref{eq:am-asymptotic} we recall that, for any numbers $\alpha>\beta>0$ it holds that
    \[
    \frac{\Gamma(m+\alpha)}{\Gamma(m+\beta)}=m^{\alpha-\beta}(1+O(m^{-1})),
    \]
    hence, recalling also the definition of the Pochhammer symbol, we find 
    \[
    a_m^{(n,s)}=|\mathbb S^{n-1}|\frac{ \Gamma(\frac n2)\Gamma(m+\frac{n+s}{2})\Gamma(m+1+\frac s2)}{\Gamma(\frac{n+s}{2})\Gamma(1+\frac s2)\Gamma(m+\frac n2)\Gamma(m+1)}=\An m^s+ O(m^{s-1})
    \]
    where $\An$ is the one defined in \eqref{eq:am-asymptotic}.
\end{proof}
Thanks to Lemma \ref{lem:angular} we can derive a power series expansion for the interaction kernels $I$ and $J.$\begin{lemma}
    \label{espansione in serie I e J}
    For any $0<\rho<r<R$ it holds that \begin{equation}\label{asintotico J}
    J(r,R)=\sum_{m=0}^{+\infty}b_mr^{2m}R^{-2m-s}\quad \text{with}\quad b_m=\frac{a_m^{(n,s)}}{2m+s}
    ,\end{equation}and\begin{equation}\label{asintotico I}
     I(r,\rho)=\sum_{m=0}^{+\infty}c_m\rho^{n+2m}r^{-n-s-2m}\quad \text{with}\quad c_m=\frac{a_m^{(n,s)}}{2m+n}.
    \end{equation}In particular there exists constants $0<c_{n,s}<C_{n,s}$ that depends only on $n$ and $s$ such that \begin{equation} \label{asintotici b e c}
    c_{n,s}(m+1)^{s-1}\leq b_m\leq C_{n,s}(m+1)^{s-1},\quad c_{n,s}(m+1)^{s-1}\leq c_m\leq C_{n,s}(m+1)^{s-1},
    \end{equation}and \begin{equation}\label{asintotico b-c}
     |b_m-c_m|\le C_{n,s}(m+1)^{s-2}.
    \end{equation}
\end{lemma}
\begin{proof}
    We prove only the expansion of $J,$ with the very same computation we obtain the analogues result for $I.$
    Let us write $x_r=re_1,$ by polar coordinates we obtain \[\begin{aligned}
    J(r,R)=\int_{\R^n\setminus B_R}\frac{dy}{|x_r-y|^{n+s}}&=\int_{R}^{+\infty}\int_{\mathbb S^{n-1}}\frac{\rho^{n-1}}{|re_1-\rho\omega|^{n+s}}\dd\omega\dd\rho\\&=\int_{R}^{+\infty}\frac{\dd\rho}{\rho^{1+s}}\int_{\mathbb S^{n-1}}\frac{\dd\omega}{(1-2\frac{r}{\rho}\omega_1+\frac{r^2}{\rho^2})^\frac{n+s}{2}},
  \end{aligned}  \]where the last equality simply follows by algebraic manipulations and by observing that \[
    |x_r-y|^2=|re_1-\rho\omega|^2=r^2+\rho^2-2r\rho\omega_1.
    \]Thanks to the definition of $\mathscr A_{n,s}$ and Lemma \ref{lem:angular} we can rewrite this last integral as 
\[
J(r,R)=\int_R^{+\infty}\rho^{-1-s}\mathscr A_{n,s}\bigg(\frac{r}{\rho}\bigg)\dd\rho=\int_{R}^{+\infty}\rho^{-1-s}\sum_{m=0}^{+\infty}a_m^{(n,s)}\left(\frac{r}{\rho}\right)^{2m}.
\]Since the coefficients are positive, we can interchange the order of the series and the integral, obtaining\[
J(r,R)=\sum_{m=0}^{+\infty}a_m^{(n,s)}r^{2m}\int_R^{+\infty}\rho^{-1-s-2m}\dd\rho=\sum_{m=0}^{+\infty}\frac{a_m^{(n,s)}}{2m+s}r^{2m}R^{-2m-s}
\]that is exactly \eqref{asintotico J}. With similar computations we also obtain the expansion of \eqref{asintotico I}. Finally the asymptotic behaviors in \eqref{asintotici b e c} and \eqref{asintotico b-c} simply follow from \eqref{eq:am-asymptotic}.
\end{proof}
The following elementary results will be needed in the asymptotic analysis of the fractional curvature.
\begin{lemma}The following facts hold true:\label{lem:laplace}

 \noindent Let $\sigma>-1$ and $c>0$, then for all $0<h\le1$,
\begin{equation}\label{eq:laplace-upper}
 \sum_{m=0}^{+\infty}(m+1)^\sigma e^{-c(m+1)h}
 \le C_{\sigma,c}h^{-1-\sigma}.
\end{equation}
 For any $\tau>1$ and $K\ge1$ we have,
\begin{equation}\label{eq:power-tail}
 \sum_{m= K}^{+\infty}(m+1)^{-\tau}\le C_\tau K^{1-\tau}.
\end{equation}
Let $A_m$ be a sequence such that 
\begin{equation}\label{assunzione}A_m=A(m+1)^\sigma+O((m+1)^{\sigma-1})\quad\text{for some $\sigma>-1$}.\end{equation} Then for any $\varphi\in C^0((0,\infty))$ such that
\begin{equation}\label{assunzioni}
 |u^\sigma\varphi(u)|\le
 C\begin{cases}
 u^{\eta_0},&0<u\le1,\\
 u^{-1-\eta_\infty},&u\ge1,
 \end{cases}
\end{equation}for some
$\eta_0>-1$, $\eta_\infty>0$, and $C>0$ it holds that
\begin{equation}\label{eq:riemann-general}
\lim_{h\to0^+}\: h^{1+\sigma}\sum_{m=0}^{+\infty} A_m\varphi((m+1)h)
= A\int_0^{+\infty} u^\sigma\varphi(u)\,\dd u.
\end{equation}
\end{lemma}
\begin{proof}
    The estimates in \eqref{eq:laplace-upper} and \eqref{eq:power-tail} easily follow by comparing each sum with the corresponding integral. We now focus on \eqref{eq:riemann-general}. Let $F(u)=Au^\sigma \varphi(u),$ by assumption \eqref{assunzioni} $F\in L^1(0,+\infty).$ For any $h>0$ we define $F_h(u)=F(h\lceil \frac{u}{h}\rceil)$, since $F$ is continuous we have that $F_h(u)\to F(u)$ as $h\to0$ for any positive $u$. Also, on every interval of the kind $(mh,(m+1)h]$ $F_h(u)$ is equal to $F((m+1)h)$. Hence we have \[
    \int_0^{+\infty} F_h(u)\dd u=h\sum_{m=0}^{+\infty}F((m+1)h)=h^{1+\sigma}\sum_{m=0}^{+\infty}A(m+1)^\sigma\varphi((m+1)h).
    \]Now we observe that $|F_h(u)|\leq C\min\{ u^{\min\{\eta_0,0\}},u^{-1-\eta_\infty}\} $ and this last function is in $L^1(0,+\infty)$. Hence from the dominated convergence theorem we obtain \[
    \int_0^{+\infty}F(u) \dd u=\lim_{h\to0^+}\int_0^{+\infty}F_h(u)\dd u=\lim_{h\to0^+}h^{1+\sigma}\sum_{m=0}^{+\infty}A(m+1)^\sigma\varphi((m+1)h).
    \]
    So we need only to show that \[
    \lim_{h\to 0^+}h^{1+\sigma}\sum_{m=0}^{+\infty}O((m+1)^{\sigma-1})\varphi((m+1)h)=0.
    \]In order to do this we split the sum where $(m+1)h\leq 1$ and where $(m+1)h\geq1 $ to obtain \[\begin{aligned}
\bigg|h^{1+\sigma}\sum_{m=0}^{+\infty}O((m+1)^{\sigma-1})\varphi((m+1)h)\bigg|&\leq Ch^{1+\sigma}\sum_{m=0}^{+\infty}(m+1)^{\sigma-1}|\varphi((m+1)h)|\\
&\leq Ch^{1+\eta_0}\sum_{(m+1)h\leq 1}(m+1)^{\eta_0-1}\\&\hspace{1em}+Ch^{-\eta_\infty}\sum_{(m+1)h\geq 1}(m+1)^{-2-\eta_\infty}.
    \end{aligned}\]For the sum where $(m+1)h\leq 1$ we distinguish three cases: if $\eta_0<0$ the sum converges and hence the limit as $h$ goes to $0$ vanishes. If $\eta_0=0$ by comparison with the integral we find \[
    Ch^{1+\eta_0}\sum_{(m+1)h\leq 1}(m+1)^{\eta_0-1}\leq C h^{1+\eta_0}\log\bigg(\frac{1}{h}\bigg)\to0,
    \]while if $\eta_0>0$ we simply observe that \[
   h^{1+\eta_0} \sum_{(m+1)h\leq 1}(m+1)^{\eta_0-1}\leq C h\to 0.
    \]In the same way we estimate \[
    h^{-\eta_\infty}\sum_{(m+1)h\geq 1}(m+1)^{-2-\eta_\infty}\leq Ch\to0,
    \]that concludes the proof.
\end{proof}
In order to derive an asymptotic expansion for $k_{n,s}$ we take advantage of the following auxiliary function.
For $q\in(0,1)$ define
\begin{equation}\label{eq:Phi}
 \Phi_{n,s}(q):=\sum_{m=0}^{+\infty}\left[
 b_m\frac{q^{2m+s}}{1+q^{2m+s}}
 -c_m\frac{q^{2m+n}}{1+q^{2m+n}}
 \right].
\end{equation}

\begin{lemma}\label{lem:Phi}
Every summand in \eqref{eq:Phi} is positive and increasing in $q$.  Moreover,
\begin{equation}\label{eq:Phi-limit}
 \Phi_{n,s}(1^-):=\lim_{q\to1^-}\Phi_{n,s}(q)=\frac12\sum_{m=0}^{+\infty}(b_m-c_m)<+\infty.
\end{equation}
\end{lemma}
\begin{proof}
    For any $q\in(0,1)$ we define the function $f_q:(0,+\infty)\to\R$ as \[
    f_q(t)=\frac{1}{t}\frac{q^t}{1+q^t},
    \]and we observe that \[
    b_m\frac{q^{2m+s}}{1+q^{2m+s}}- c_m\frac{q^{2m+n}}{1+q^{2m+n}}=a_m^{(n,s)}(f_q(2m+s)-f_q(2m+n)).
    \]So in order to prove that every summand is positive it is enough to show that $f_q$ is strictly decreasing, but this just follows by observing that \[
    f'_q(t)=-\frac{1}{t^2}\frac{q^t}{1+q^t}+\frac{1}{t}\frac{q^t\log q}{(1+q^t)^2}<0.
    \]Now we define, \[
    T_m(q)=b_m\frac{q^{2m+s}}{1+q^{2m+s}}-c_m\frac{q^{2m+n}}{1+q^{2m+n}}
    \]and we compute \[
    T_m'(q)=a_m^{(n,s)}\left[\frac{q^{2m+s-1}}{(1+q^{2m+s})^2}-\frac{q^{2m+n-1}}{(1+q^{2m+n})^2}\right]=\frac{a_m^{(n,s)}}{q}[\psi(q^{2m+s})-\psi(q^{2m+n})]
    \]where $\psi(z)={z}{(1+z)^{-2}}.$ It is easy to see that $\psi$ is increasing and hence we have $T_m'(q)>0$ and hence every summand is increasing in $q.$ From thees observations we obtain that $
    T_m(q)\leq \frac12(b_m-c_m),
    $ that is summable by \eqref{asintotico b-c}. Therefore from the dominated convergence theorem we infer \[
    \lim_{q\to1^-}\Phi_{n,s}(q)=\frac12\sum_{m=0}^{+\infty}(b_m-c_m).
    \]
\end{proof}
The next proposition shows that the limit of the function $\Phi$ as $q\to 1^-$ is half of the fractional curvature.
\begin{proposition}\label{prop:critical-identity}
We have 
\begin{equation}\label{eq:critical-identity}
 2\Phi_{n,s}(1^-)=\sum_{m=0}^{+\infty}(b_m-c_m)=\kappaN
\end{equation}
\end{proposition}
\begin{proof}
    The first equality was shown in Lemma \ref{lem:Phi} so we need only to prove the second one. From \eqref{asintotico J} and \eqref{asintotico I} we have \[
    b_m-c_m=\frac{(n-s)a_m}{(2m+s)(2m+n)}=\frac{(n-s)|\mathbb S^{n-1}|(\frac{n+s}{2})_m(\frac{s}{2})_m}{sn(\frac{n}{2}+1)_mm!}
    \]where the last equality follows from the definition of $a_m$ in \eqref{eq:angular-series} and by observing that \[
    \frac{1}{2m+s}=\frac{(\frac{s}{2})_m}{s(1+\frac s2)_m}\quad\text{and}\quad \frac{1}{2m+n}=\frac{(\frac{n}{2})_m}{n(1+\frac n2)_m}.
    \]By the definition of hypergeometric function and \cite[Section 1.3]{Bai64} we find \begin{equation}\begin{aligned}
    \sum_{m=0}^{+\infty}(b_m-c_m)=\frac{(n-s)|\mathbb S^{n-1}|}{sn}\sum_{m=0}^{+\infty}\frac{(\frac{n+s}{2})_m(\frac{s}{2})_m}{(\frac{n}{2}+1)_mm!}&=\frac{(n-s)|\mathbb S^{n-1}|}{sn}{}_2F_1\big(\frac{n+s}{2},\frac s2,\frac n2+1,1\big)\\
    &=\frac{(n-s)|\mathbb S^{n-1}|}{sn}\frac{\Gamma(\frac n2+1)\Gamma(1-s)}{\Gamma(1-\frac s2)\Gamma(\frac{n-s}{2}+1)}\\
    &=\frac{2\pi^{\frac n2}\Gamma(1-s)}{s\Gamma(1-\frac s2)\Gamma(\frac{n-s}{2})}\label{ultima}.
   \end{aligned} \end{equation}From the Legendre duplication formula for the $\Gamma$ function (see \cite[Theorem 1.5.1]{AAR}) we have \[
   \sqrt\pi\Gamma(1-s)=2^{-s}\Gamma\big(\frac{1-s}{2}\big)\Gamma\big(1-\frac s2\big)
   \]that combined with \eqref{ultima} and Lemma \ref{k esplicito} gives \[
   \sum_{m=0}^{+\infty}(b_m-c_m)=\frac{2\pi^{\frac n2}\Gamma(1-s)}{s\Gamma(1-\frac s2)\Gamma(\frac{n-s}{2})}
 =\frac{2^{1-s}\pi^{(n-1)/2}\Gamma(\frac{1-s}{2})}
 {s\,\Gamma(\frac{n-s}{2})}=\kappaN.
   \]
\end{proof}
Now we derive a first order expansion of the function $\Phi_{n,s}$ in a left neighborhood of $1.$ For the sake of readability we introduce the following notations. For $u\geq0$ we define the functions
\begin{equation}\label{eq:fgh}
 f(u):=\frac{e^{-u}}{1+e^{-u}},\qquad
 g(u):=\frac{e^{-u}}{(1+e^{-u})^2},\qquad
 h(u):=\frac{e^{-u}(1-e^{-u})}{(1+e^{-u})^3},
\end{equation}
and the constant
\begin{equation}\label{eq:Dstar}
 D_{n,s,*}:=\frac{\An}{2^{1+s}}\int_0^{+\infty} u^sh(u)\,\dd u>0.
\end{equation}where $\An$ is the one in \eqref{eq:am-asymptotic}.
We are finally ready to derive the asymptotic expansion of the fractional mean curvature.
\begin{proposition}\label{prop:critical-defect}
As $p\to0^+$ the following expansion holds,
\begin{equation}\label{eq:critical-defect}
 \Phi_{n,s}(1^-)-\Phi_{n,s}(e^{-p})
 =\frac{n-s}{1-s}D_{n,s,*}p^{1-s}+o(p^{1-s}).
\end{equation}
In particular, there are $c_{n,s}<C_{n,s}$ and $p_0>0$ such that
\begin{equation}\label{eq:defect-two-sided}
 c_{n,s}p^{1-s}\le\Phi_{n,s}(1^-)-\Phi_{n,s}(e^{-p})
 \le C_{n,s}p^{1-s},
 \qquad 0<p<p_0.
\end{equation}
\end{proposition}
\begin{proof}
    By Lemma \ref{lem:Phi} we have \[\begin{aligned}
    \Phi_{n,s}(1^-)-\Phi_{n,s}(e^{-p})&=\frac{1}{2}\sum_{m=0}^{+\infty}(b_m-c_m)-\sum_{m=0}^{+\infty}\left[b_m \frac{e^{-p(2m+s)}}{1+e^{-p(2m+s)}}-c_m\frac{e^{-p(2m+n)}}{1+e^{-p(2m+n)}}\right]\\
    &=\frac{1}{2}\sum_{m=0}^{+\infty}(b_m-c_m)-\sum_{m=0}^{+\infty}\left[b_m f(p(2m+s))-c_mf(p(2m+n))\right].
   \end{aligned} \]If we define the function $G:[0,+\infty)\to\R$  as\[
   G(u)=\frac{\frac12-f(u)}{u},\quad G(0)=\frac{1}{4},
   \]we can rewrite the above relation as \[
     \Phi_{n,s}(1^-)-\Phi_{n,s}(e^{-p})
 =p\sum_{m=0}^{+\infty} a_m^{(n,s)}
 \bigl[G(p(2m+s))-G(p(2m+n))\bigr].
   \]It is easy to see that the function $G$ belongs to $C^2([0,+\infty))$ with second derivative that decay like $u^{-3}$ as $u\to+\infty$. Hence, by a Taylor expansion in the point $2mp$ we find\[
G(p(2m+s))=G(2mp)+spG'(2mp)+p^2G''(\xi_1)
   \]for some $\xi_1\in (2mp,p(2m+s))$, and\[
   G(p(2m+n))=G(2mp)+npG'(2mp)+p^2G''(\xi_2)
   \]for some $\xi_2\in (2mp,p(2m+n)).$ That together implies
  \[  G(p(2m+s))-G(p(2m+n))
 =-p(n-s)G'(2mp)+\mathcal R_{m,p},
\]
with the remainder that can be bounded as
\[
 |\mathcal R_{m,p}|\le C_{n,s}
 \begin{cases}
  p^2&\quad\text{if }mp\le1,\\
  \frac{1}{pm^3}
 &\quad\text{if }mp>1.
 \end{cases}
\]
By \eqref{eq:am-asymptotic} we can infer that $a_m^{(n,s)}\leq C_{n,s}(m+1)^s,$ hence we have \[\begin{aligned}
p\sum_{m=0}^{+\infty}a_m^{(n,s)}|\mathcal R_{m,p}|&=p\sum_{mp\leq1}a_m^{(n,s)}|\mathcal R_{m,p}|+p\sum_{mp\geq1}a_m^{(n,s)}|\mathcal R_{m,p}|\\
&\leq C_{n,s}p^3\sum_{mp\leq1}(m+1)^s+C_{n,s}\sum_{mp\geq1}m^{s-3}\\
&\leq C_{n,s}p^{2-s}=o(p^{1-s})\quad\text{as $p\to0$.}\end{aligned}
\]
Therefore
\[
 \Phi_{n,s}(1^-)-\Phi_{n,s}(e^{-p})
 =-(n-s)p^2\sum_{m=1}^{+\infty} a_m^{(n,s)}G'(2mp)+o(p^{1-s}),
\]note that we can neglect the first term of the summation since, as one can easily compute, $G'(0)=0.$ By \eqref{eq:am-asymptotic} and some standard computations, it is easy to see that the assumptions \eqref{assunzione} and \eqref{assunzioni} in Lemma \ref{lem:laplace} are fulfilled, hence we obtain
\[
 p^{1+s}\sum_{m=0}^{+\infty} a_m^{(n,s)}G'(2mp)
 \longrightarrow \frac{\An}{2^{1+s}}
 \int_0^{+\infty} u^sG'(u)\,\dd u,
\]that implies \[
\Phi_{n,s}(1^-)-\Phi_{n,s}(e^{-p})
 =-(n-s)\frac{\An}{2^{1+s}}p^{1-s}\int_0^{+\infty}u^sG'(u)\dd u+o(p^{1-s}).
\]By integrating by parts and recalling that $uG(u)=\frac{1}{2}-f(u)$ we obtain \[\begin{aligned}
-\int_0^{+\infty}u^s G'(u)\dd u&=s\int_0^{+\infty}u^{s-1}G(u)\dd u=s\int_0^{+\infty}u^{s-2}\left(\frac{1}{2}-f(u)\right)\dd u\\
&=s\int_0^{+\infty}u^{s-2}\left(\int_0^ug(t)\dd t\right) du
\end{aligned}\]where the last equality follows by observing that $f'(u)=-g(u).$
By Fubini's theorem,\[
s\int_0^{+\infty}u^{s-2}\left(\int_0^ug(t)\dd t\right)=s\int_0^{+\infty}g(t)\left(\int_t^{+\infty}u^{s-2}\dd u\right)\dd t=\frac{s}{1-s}\int_0^{+\infty}t^{s-1}g(t)\dd t,
\]hence \[
-\int_0^{+\infty}u^sG'(u)\dd u=\frac{s}{1-s}\int_0^{+\infty}u^{s-1}g(u)\dd u.
\]Since $g'(u)=-h(u)$ another integration by parts shows that\[
{s}\int_0^{+\infty}u^{s-1}g(u)\dd u=-\int_0^{+\infty}u^sg'(u)\dd u=\int_0^{+\infty}u^s h(u)\dd u,
\]that concludes the proof of \eqref{eq:critical-defect}. Finally, \eqref{eq:defect-two-sided} immediately follows from \eqref{eq:critical-defect}.
\end{proof}
\subsection{Geometric separation properties}\label{sec:separation}
The aim of this subsection is to prove that properties \eqref{eq:central-main} and \eqref{eq:adjacent-main} of Theorem \ref{thm:main finite} hold true. We begin with a simple Lemma.
\begin{lemma}\label{lem:bathtub}
Let $0<\sigma<1$, $d>0$, and let $E\subset(0,+\infty)$ be measurable.  If
\begin{equation}\label{integrale}
 \int_E \sigma\,t^{-1-\sigma}\,\dd t<+\infty,
\end{equation}
then
\begin{equation}\label{eq:bathtub}
 \int_E \sigma\,(t+d)^{-1-\sigma}\,\dd t
 \le\bigg(d+\left(\int_E\sigma\, t^{-1-\sigma}\dd t\right) ^{-\frac1\sigma}\bigg)^{-\sigma}.
\end{equation}
\end{lemma}
\begin{proof}
    If the integral in \eqref{integrale} is $0$, the conclusion is immediate. Then we will assume that the integral in \eqref{integrale} is strictly positive. Let $\nu$ be the measure whose density (with respect to the Lebesgue measure on $(0,+\infty)$) is $\sigma t^{-1-\sigma}.$ For any $T>0$ we have \[
    \nu([T,+\infty))=T^{-\sigma},
    \]hence, if we choose $T$ such that $T^{-\sigma}=\int_E \sigma t^{-1-\sigma}\dd t$ we find \[
     \nu([T,+\infty))=\int_ E \sigma t^{-1-\sigma}\dd t=\nu(E).
    \]Moreover, if we define $w_d(t)=\left(\frac{t}{t+d}\right)^{1+\sigma}$, we have \[
    \int_E \sigma(t+d)^{-1-\sigma}\,dt=\int_E w_d(t)\dd\nu(t).
    \]Now let $E^-=E\cap (0,T)$ and $E^+=[T,+\infty)\setminus E$, since $\nu(E)=\nu([T,+\infty))$ we easily have that $\nu(E^-)=\nu(E^+)$. By monotonicity we obtain \[
    \int_{E^-}w_d(t)\dd\nu(t)\leq w_d(T)\nu(E^-)\quad\text{and}\quad \int_{E^+}w_d(t)\dd\nu(t)\geq w_d(T)\nu(E^+),
    \]that together with the fact that $\nu(E^-)=\nu(E^+)$ implies \[
    \int_{E^-}w_d(t)\dd\nu(t)\leq \int_{E^+}w_d(t)\dd\nu(t).
    \]Now we observe that \[
    \begin{aligned}
    \int_E\sigma(t+d)^{-1-\sigma}\dd t=\int_E w_d(t)\dd\nu(t)&=\int_{E^-}w_d(t)\dd\nu(t)+\int_{E\cap[T,+\infty)}w_d(t)\dd\nu(t)\\
    &\leq\int_{E^+}w_d(t)\dd\nu(t)+\int_{E\cap[T,+\infty)}w_d(t)\dd\nu(t)\\
    &=\int_T^{+\infty}w_d(t)\dd\nu(t)=\int_T^{+\infty}\sigma(t+d)^{-1-\sigma}\dd t =(T+d)^{-\sigma}
    ,\end{aligned}
    \]the conclusion follows by recalling the choice of $T.$
\end{proof}
\begin{lemma}\label{lem:screening}
Let $r_k>0$, $d_k>0$ such that $d_k/r_k\to0$.  Let
$E_k\subset(r_k+d_k,+\infty)$ be measurable and set
\[
 M_k^+:=\int_{E_k}\mathcal K(r_k+d_k,\rho)\,\dd\rho,
 \qquad
 M_k^0:=\int_{E_k}\mathcal K(r_k,\rho)\,\dd\rho.
\]
If
\[
 M_k^+=\GammaN d_k^{-s}+o(d_k^{-s}),
\]
then
\begin{equation}\label{eq:screening-out}
 M_k^0\le(2^{-s}+o(1))\GammaN d_k^{-s}.
\end{equation}The same result holds also if $E_k\subset(0,r_k)$ by changing the role of $M_k^+$ and $M_k^0,$ namely, if  \[
 M_k^0=\GammaN d_k^{-s}+o(d_k^{-s}),
\]then \begin{equation}\label{eq:screening-in}
M_k^+\le(2^{-s}+o(1))\GammaN d_k^{-s}.
\end{equation}
\end{lemma}
\begin{proof}
   We only prove \eqref{eq:screening-out}, the proof for \eqref{eq:screening-in} is similar. Let $\delta\in(0,\frac{\delta_0}{4})$ where $\delta_0$ is the one of Lemma \ref{lem:near-contact} and $\widetilde E_k=\{t>0\,:\: r_k+d_k+t\in E_k\}$. By a change of variable we can write \[
    M^+_k=\int_{\widetilde E_k}  \mathcal K(r_k+d_k,r_k+d_k+t)\dd t\quad\text{and}\quad M^0_k=\int_{\widetilde E_k}  \mathcal K(r_k,r_k+d_k+t)\dd t.
    \]Since $\frac{d_k}{r_k}\to 0$ we can assume $d_k\leq \delta r_k$. Let $\widetilde E^\delta_k=\widetilde E_k\cap(0,\delta r_k)$, by \eqref{eq:Kglobal} we have \begin{equation}\label{mk}
    \int_{\widetilde E_k\setminus \widetilde E_k^\delta} \mathcal K(r_k+d_k,r_k+d_k+t)\dd t\leq C_{n,s}\int_{\delta r_k}^{+\infty}t^{-1-s}\dd t= C_{n,s,\delta}r_k^{-s}=o(d_k^{-s}),
    \end{equation}and in the same way \begin{equation}\label{m0}
     \int_{\widetilde E_k\setminus \widetilde E_k^\delta} \mathcal K(r_k,r_k+d_k+t)\dd t\leq C_{n,s}\int_{\delta r_k}^{+\infty}(d_k+t)^{-1-s}\dd t\leq C_{n,s,\delta}r_k^{-s}=o(d_k^{-s}).
    \end{equation}Since $d_k=o(r_k)$ and by \eqref{eq:density-expansion} we find \[
    \mathcal K(r_k+d_k,r_k+d_k+t)=\gamma_{n,s}t^{-1-s}+O(r_k^{-1}t^{-s})
    \]and \[
    \mathcal K(r_k,r_k+d_k+t)=\gamma_{n,s}(t+d_k)^{-1-s}+O(r_k^{-1}(t+d_k)^{-s}).
    \]A straightforward computation shows that \[
    r_k^{-1}\int_0^{\delta r_k}t^{-s}\dd t=o(d_k^{-s})\quad\text{and}\quad r_k^{-1}\int_0^{\delta r_k}(t+d_k)^{-s}\,dt=o(d_k^{-s}),
    \]hence, taking also into account \eqref{mk} and \eqref{m0} and recalling that $\gamma_{n,s}=s\Gamma_{n,s}$ we find \[
    M_k^+=\Gamma_{n,s}\int_{\widetilde E_k^\delta}st^{-1-s}\dd t+o(d_k^{-s})\quad\text{and}\quad M_k^0=\Gamma_{n,s}\int_{\widetilde E_k^{\delta} }s(t+d_k)^{-1-s}\dd t+o(d_k^{-s}).
    \]By Lemma \ref{lem:bathtub} and the assumption on $M_k^+$ we can finally conclude that \[
    \begin{aligned}
        M_k^0=\Gamma_{n,s}\int_{\widetilde E_k^\delta}s(t+d_k)^{-1-s}\dd t+o(d_k^{-s})&\leq \Gamma_{n,s}\bigg(d_k+\bigg(\int_{\widetilde E_k^\delta}st^{-1-s}\dd t\bigg)^{-\frac{1}{s}}\bigg)^{-s}\\
        &=\Gamma_{n,s}(2d_k+o(d_k))^{-s}\\
        &=\Gamma_{n,s}d_k^{-s}(2^{-s}+o(1)),
    \end{aligned}
    \]that is \eqref{eq:screening-out}.
\end{proof}
We are finally ready to prove the first part of Theorem \ref{thm:main finite}, namely equation \eqref{eq:central-main}.
\begin{theorem}\label{thm:central}
There exists $c_{n,s}>0$ such that every finite radial equilibrium satisfies
\begin{equation}\label{raggi lontani}
 r_2-r_1\ge c_{n,s}r_1.
\end{equation}
\end{theorem}
\begin{proof}
    We start with the case $M=2.$ In this case the balance equation \eqref{eq:balanced} reduces to \[J(r_1,r_2)=\frac{\kappa_{n,s}}{2r_1^s}.\]By defining $p_1=\log \frac{r_2}{r_1}$ and recalling that $J(r_1,r_2)=r_1^{-s} \mathcal J(\frac{r_2}{r_1})$ we can rewrite the above equation as \[
    \mathcal J(e^{p_1})=\frac{\kappa_{n,s}}{2}.
    \]As a function of $p$, 
$\mathcal J(e^p)$ is strictly decreasing, goes to $+\infty$ as
$p\to0^+$ by Lemma \ref{lem:near-contact}, and to zero as
$p\to+\infty$. Hence $p_1$ is the unique solution to \[
\mathcal J(e^p)=\frac{\kappa_{n,s}}{2}
,\]and therefore it depends only on $n$ and $s$. This shows that \[
\frac{r_2-r_1}{r_1}=\frac{r_2}{r_1}-1=e^{p_1}-1>0
,\]that is \eqref{raggi lontani}. Let us deal with the case $M\geq 4,$ if \eqref{raggi lontani} is not true there would exist a sequence of equilibria (possibly with different numbers of interfaces) such that $r_2-r_1=o(r_1).$ If we define
\[
 \mathcal A:=\bigcup_{j=1}^{N-1}(r_{2j+1},r_{2j+2}).
\]we can observe that, by \eqref{eq:IJ}, \eqref{A+i} and \eqref{derivata}, it holds that \[
A_2^+=\int_\mathcal A \mathcal K(r_2,\rho)\dd\rho
\]while $A_2^-=I(r_2,r_1).$ By Lemma \ref{lem:near-contact} we may infer \[
I(r_2,r_1)=\Gamma_{n,s}(r_2-r_1)^{-s}+O(r_2^{-s})=\Gamma_{n,s}(r_2-r_1)^{-s}+o((r_2-r_1)^{-s})
\]where the last equality follows from the fact that $r_2-r_1=o(r_1).$ By these observations and the balance equation \eqref{eq:balanced} we find \[
A_2^+=\Gamma_{n,s}(r_2-r_1)^{-s}+o(r_2-r_1)^{-s}
\]If we define \[
 B_1:=\int_{\mathcal A}\mathcal K(r_1,\rho)\,\dd\rho
\]we can apply Lemma \ref{lem:screening} to obtain \begin{equation}
B_1 \le(2^{-s}+o(1))\GammaN (r_2-r_1)^{-s}.\label{B}
\end{equation}Now we observe that \[
B_1=\sum_{k=1}^{M-2}(-1)^{k+1}J(r_1,r_{2+k})
,\] hence from \eqref{A+i} we have \[
A_1^+=J(r_1,r_2)-B_1.
\]By Lemma \ref{lem:near-contact} and \eqref{B} we find \[\begin{aligned}
A_1^+&\geq\Gamma_{n,s}(r_2-r_1)^{-s}+o((r_2-r_1)^{-s})-(2^{-s}+o(1))\GammaN (r_2-r_1)^{-s}\\&=\Gamma_{n,s}(1-2^{-s})(r_2-r_1)^{-s}+o((r_2-r_1)^{-s}).
\end{aligned}\]That is not possible, indeed the balance equation \eqref{eq:balanced} for $i=1$ combined with the assumption $r_2-r_1=o(r_1)$ implies $A_1^+=\frac{\kappa_{n,s}}{2}r_1^{-s}=o((r_2-r_1)^{-s}).$
\end{proof}
Now we can prove the adjacent gap estimate, that is \eqref{eq:adjacent-main} in Theorem \ref{thm:main finite}.
\begin{theorem}\label{thm:adjacent}
There are $c_0,C_0,\eta_0>0$, depending only on $n$ and $s$, such that
\eqref{eq:adjacent-main} holds.
\end{theorem}
\begin{proof}
    We argue by contradiction. If the statement is false we can find a sequence of finite radial equilibria (that we do not relabel) such that $\min\{p_i,p_{i+1}\}\to 0$ and one of the ratios $p_{i+1}/p_i$ or $p_i/p_{i+1}$ goes to $0.$ If we recall that $p_i=\log\frac{r_{i+1}}{r_i}$ and define $g_i=r_{i+1}-r_i$ we can see that $g_i=r_i(e^{p_i}-1),$ hence \[
    \frac{g_{i+1}}{g_i}=\frac{e^{p_i}(e^{p_{i+1}}-1)}{e^{p_i}-1}.
    \]This relation, combined with the fact that $p_{i+1}\to0$ if ${p_{i+1}}/{p_i}\to 0$ shows that ${g_{i+1}}/{g_i}\to0$ if ${p_{i+1}}/{p_i}\to0.$ Analogously we can see that ${g_{i}}/{g_{i+1}}\to0$ if ${p_{i}}/{p_{i+1}}\to0.$

    Let us start with the case ${p_{i+1}}/{p_i}\to0,$ so ${g_{i+1}}/{g_i}\to0.$ By the fact that the summation in \eqref{A-i} is alternate we can see that \begin{equation}\label{alternata}
    I(r_{i+2},r_{i+1})-I(r_{i+2},r_i)\leq A^-_{i+2}\leq I(r_{i+2},r_{i+1}).
    \end{equation}
    Since \[
    \frac{g_{i+1}}{r_{i+1}}\leq \frac{g_{i+1}}{r_{i+1}-r_i}=\frac{g_{i+1}}{g_i}\to 0
    \]we have $g_{i+1}\to 0$ hence,
    by Lemma \ref{lem:near-contact} we have \begin{equation}\label{espansione}
    I(r_{i+2},r_{i+1})=\Gamma_{n,s}g_{i+1}^{-s}+o(g_{i+1}^{-s}).
    \end{equation}Furthermore, by \eqref{eq:Kglobal} we find \begin{equation}\label{i+2i}
    I(r_{i+2},r_i)\leq C_{n,s}(g_i+g_{i+1})^{-s}=o(g_{i+1}^{-s})
    \end{equation}where the last equality follows simply observing that \[
    \frac{(g_i+g_{i+1})^{-s}}{g_{i+1}^{-s}}=\left(\frac{g_{i+1}}{g_i+g_{i+1}}\right)^s=\left(\frac{g_{i+1}/g_i}{1+g_{i+1}/g_i}\right)^s\to 0.
    \]From \eqref{alternata}, \eqref{espansione} and \eqref{i+2i} we infer \[
    A_{i+2}^-=\Gamma_{n,s}g_{i+1}^{-s}+o(g_{i+1}^{-s}).
    \]By the balance equation \eqref{eq:balanced} and the fact that $r_{i+2}^{-s}\leq r_{i+1}^{-s}=o(g_{i+1})^{-s}$ we obtain \[
    A^+_{i+2}-A^-_{i+2}=\frac{\kappaN}{2 r_{i+2}^s}=o(g_{i+1}^{-s}),
    \]that implies \[
     A_{i+2}^+=\Gamma_{n,s}g_{i+1}^{-s}+o(g_{i+1}^{-s}).
    \]Now let \[
    B_{i+1}=\sum_{k=1}^{M-i-2}(-1)^{k+1}J(r_{i+1},r_{i+2+k}),
    \]applying Lemma \ref{lem:screening} as in the proof of Theorem \ref{thm:central} we infer \[
    B_{i+1}\leq(2^{-s}+o(1))\Gamma_{n,s}g_{i+1}^{-s},
    \]that combined with Lemma \ref{lem:near-contact} and observing that \[
    A_{i+1}^+=J(r_{i+1},r_{i+2})-B_{i+1}
    ,\]implies\begin{equation}\label{contraddizione}
    A_{i+1}^+\geq(1-2^{-s})\Gamma_{n,s}g_{i+1}^{-s}+o(g_{i+1}^{-s}).
    \end{equation}However, by \eqref{A-i} we have $A_{i+1}^-\leq I(r_{i+1},r_i)$ that combined with the balance equation \eqref{eq:balanced} and \eqref{eq:Kglobal} gives\[
    A_{i+1}^+\leq I(r_{i+1},r_i)+\frac{\kappa_{n,s}}{2 r_{i+1}^s}\leq C_{n,s}g_i^{-s}+\frac{\kappa_{n,s}}{2 r_{i+1}^s}=o(g_{i+1}^{-s})
    \]that contradicts \eqref{contraddizione}.
    
    The case $p_i/p_{i+1}\to0$ for $i\geq 2$ can be treated in a similar way, just using \eqref{eq:screening-in} in place of \eqref{eq:screening-out}. We need to deal with the case $p_1/p_2\to 0$ but this one just follows from Theorem \ref{thm:central}. Indeed if $p_1/p_2\to0$ we would have \[
    \frac{r_2-r_1}{r_1}=e^{p_1}-1\to 0,
    \]a contradiction. 
    
    It only remains to show that the constants that we found can be chosen independently of the configuration. If this is not the case, for any $k\in\mathbb N$ we can find a configuration and an index $i_k$ such that $$\min\{p_{i_k},p_{i_k+1}\}\leq \frac{1}{k}\quad\text{but}\quad\frac{p_{i_k}}{p_{i_k+1}}\notin\bigg[\frac{1}{k},k\bigg].$$ Hence either $\frac{p_{i_k+1}}{p_{i_k}}\to0$ or $\frac{p_{i_k}}{p_{i_k+1}}\to0,$ in both cases we find a contradiction by repeating the argument done before.
\end{proof}
\subsection{Approximation of the nonlocal interaction term}\label{sec:stress}In this subsection, we introduce some quantities associated with the interactions between the interfaces. These quantities allow us to approximate the nonlocal interaction term and to estimate the corresponding error.

For $1\leq i\le M-1$ and $\lambda>0$, we define the alternating responses to the right and to the left of $r_i$ as 
\begin{equation}\label{eq:RLdef}
R_i(\lambda):=\sum_{k=1}^{M-i}(-1)^{k+1}\left(\frac{r_i}{r_{i+k}}\right)^\lambda
\quad,\quad L_i(\lambda):=\sum_{k=1}^{i-1}(-1)^{k+1}\left(\frac{r_{i-k}}{r_i}\right)^\lambda
,\end{equation}we also set $R_M(\lambda)=0$ and $L_1(\lambda)=0$. Equivalently, in terms of the logarithmic gaps, they can also be written as \begin{equation}\label{logaritmo}
R_i(\lambda)=\sum_{k=1}^{M-i}(-1)^{k+1}
 \exp\bigg(-\lambda\sum_{q=i}^{i+k-1}p_q\bigg)\quad,\quad L_i(\lambda)=\sum_{k=1}^{i-1}(-1)^{k+1}
 \exp\bigg(-\lambda\sum_{q=i-k}^{i-1}p_q\bigg)
\end{equation}We notice that, by \eqref{asintotico I}, \eqref{asintotico J}, \eqref{A+i} and \eqref{A-i} we have \begin{equation}
r_i^sA_i^+=\sum_{m=0}^{+\infty}b_mR_i(2m+s)\quad,\quad r_i^s A_i^-=\sum_{m=0}^{+\infty}c_mL_i(2m+n).
\end{equation}These two relations, combined with the balance equation \eqref{eq:balanced}, and Proposition \ref{prop:critical-identity} give \begin{equation}\label{somma Phi}
\sum_{m=0}^{+\infty}b_mR_i(2m+s)-\sum_{m=0}^{+\infty}c_mL_i(2m+n)=\frac{\kappaN}{2}=\Phi_{n,s}({1^-}).
\end{equation}
We now define the interaction between the first \(i\) interfaces and the remaining ones by
\begin{equation}\label{eq:Sigma}
\Sigma_i(\lambda):=\sum_{a=1}^{i}\sum_{b=i+1}^{M}(-1)^{b-a+1}\left(\frac{r_a}{r_b}\right)^\lambda=
 \sum_{a=1}^{i}\sum_{b=i+1}^{M}(-1)^{b-a+1}
 \exp\left(-\lambda\sum_{q=a}^{b-1}p_q\right),
\end{equation}for $1\leq i\leq M-1$ and $\Sigma _M(\lambda)=\Sigma_0(\lambda)=0.$

The following Lemma summarizes some useful identities that involve the quantities $R_i(\lambda),$ $L_i(\lambda)$ and $\Sigma _i(\lambda).$ 
\begin{lemma}\label{lem:mode}
For every $\lambda>0$ it holds
\begin{equation}\label{eq:recursions}
R_i(\lambda)=e^{-\lambda p_i}(1-R_{i+1}(\lambda)),
 \qquad
 L_{i+1}(\lambda)=e^{-\lambda p_i}(1-L_i(\lambda)).
\end{equation}
Moreover,
\begin{equation}\label{eq:Sigma-div}
\Sigma_i(\lambda)-\Sigma_{i-1}(\lambda)=R_i(\lambda)-L_i(\lambda),
 \qquad
 \Sigma_i(\lambda)=R_i(\lambda)(1-L_i(\lambda)).
\end{equation}
In particular,
\begin{equation}\label{eq:mode-bounds}
 0\le R_i(\lambda)\leq e^{-\lambda p_i},\qquad 0\leq L_i(\lambda)\le e^{-\lambda p_{i-1}},
 \qquad
 0\le\Sigma_i(\lambda)\le e^{-\lambda p_i}.
\end{equation}
Furthermore, if we define
\begin{align}
 \mathcal M_i&:=\frac12\sum_{m=0}^{+\infty}
 \bigl(b_m\Sigma_i({2m+s})+c_m\Sigma_i({2m+n})\bigr),
 \label{eq:Mstress}\\
 \mathcal B_i&:=\frac12\sum_{m=0}^{+\infty}\left[
 b_m(R_i({2m+s})+L_i({2m+s}))
 -c_m(R_i({2m+n})+L_i({2m+n}))\right]
 \label{eq:Bphase}
\end{align}
the two series above converge absolutely. If the configuration is an equilibrium and $1\leq i\leq M-1$ we have
\begin{equation}\label{eq:stress-balance}
 \mathcal M_i-\mathcal M_{i-1}=\Phi_{n,s}(1^-)-\mathcal B_i.
\end{equation}
\end{lemma}
\begin{proof}
    The two identities in \eqref{eq:recursions} easily follow from \eqref{logaritmo} where $R_i$ and $L_i$ are expressed in terms of the logarithmic gaps. For \eqref{eq:Sigma-div} it is enough to observe that \[
    \Sigma_i(\lambda)-\Sigma_{i-1}(\lambda)=\sum_{b=i+1}^{M}(-1)^{b-i+1}\left(\frac{r_i}{r_b}\right)^\lambda-\sum_{a=1}^{i-1}(-1)^{i-a+1}\left(\frac{r_a}{r_i}\right)^\lambda=R_i(\lambda)-L_i(\lambda),
    \]that is the first identity. For the second one we rewrite $R_i(\lambda)$ as \[
    R_i(\lambda)=\sum_{b=i+1}^{M}(-1)^{b-i+1}\exp\bigg(-\lambda\sum_{q=i}^{b-1}p_q\bigg)
    \]and we notice that \[
    1-L_i(\lambda)=\sum_{a=1}^{i}(-1)^{i-a}\exp\bigg(-\lambda \sum_{q=a}^{i-1}p_q\bigg)
    \]where the sum in the exponent is understood to be zero if $a=i$. Therefore \[\begin{aligned}
    R_i(\lambda)(1-L_i(\lambda))&=\sum_{a=1}^i\sum_{b=i+1}^{M}(-1)^{b-a-1}\exp\bigg(-\lambda\sum_{q=i}^{b-1}p_q\bigg)\exp\bigg(-\lambda \sum_{q=a}^{i-1}p_q\bigg)\\&=\sum_{a=1}^i\sum_{b=i+1}^{M}(-1)^{b-a-1}\exp\bigg(-\lambda \sum_{q=a}^{b-1}p_q\bigg)=\Sigma_i(\lambda).
   \end{aligned} \]The first two relations in \eqref{eq:mode-bounds} follow from the fact that the sums defining $R_i$ and $L_i$ are alternating with positive and decreasing coefficients, the third one follows from the first two and the second in \eqref{eq:Sigma-div}. The fact that the series defining $\mathcal M_i$ and $\mathcal B_i$ are absolutely convergent follows from \eqref{eq:mode-bounds} and \eqref{asintotici b e c}. Finally, from the first relation in \eqref{eq:Sigma-div} and \eqref{somma Phi} we have \[\begin{aligned}
   \mathcal M_i-\mathcal M_{i-1}&=\frac{1}{2}\sum_{m=0}^{+\infty}b_m[R_i(2m+s)-L_i(2m+s)]+\frac{1}{2}\sum_{m=0}^{+\infty}c_m[R_i(2m+n)-L_i(2m+n)]\\
   &=\frac{\Phi_{n,s}(1^-)}{2}+\frac{1}{2}\sum_{m=0}^{+\infty}c_mR_i(2m+n)-b_mL_i(2m+s)\\
   &=\Phi_{n,s}(1^-)-\mathcal{B}_i,
 \end{aligned}  \]that is \eqref{eq:stress-balance}.
\end{proof}
We collect in the following lemma some identities and estimates that will be useful in the rest of the subsection.
\begin{lemma}\label{lem:spectral-cancel}
There exists a constant $C_{n,s}>0$ such that the following facts hold true
\begin{enumerate}[label=(\roman*)]
\item If $f\in C^1(0,+\infty)$  and there exists $A>0$ such that for any $\lambda>0$
$|f(\lambda)|+\lambda|f'(\lambda)|\le A$, then
\begin{equation}\label{eq:cancel-bounded}
 \left|\sum_{m=0}^{+\infty}\bigl[b_mf(2m+s)-c_mf(2m+n)\bigr]\right|
 \le C_{n,s}A.
\end{equation}
\item If $f_h\in C^1(0,+\infty)$ is a family of functions for which there exists $A>0$ and $c>0$ with
\[
 |f_h(\lambda)|\le A\lambda h e^{-c\lambda h},
 \qquad
 |f_h'(\lambda)|\le Ah(1+\lambda h)e^{-c\lambda h}
,\]{for any $\lambda\in(0,+\infty)$ and $0<h\leq 1$}
then
\begin{equation}\label{eq:cancel-small}
 \left|\sum_{m=0}^{+\infty}\bigl[b_mf_h(2m+s)-c_mf_h(2m+n)\bigr]\right|
 \le C_{n,s,c}Ah^{1-s}.
\end{equation}
\end{enumerate}
\end{lemma}
\begin{proof}
    Let us start with \eqref{eq:cancel-bounded}, for any $m\geq 0$ we have \[
    \begin{aligned}
        |b_mf(2m+s)-c_mf(2m+n)|&\leq|(b_m-c_m)f(2m+s)|+c_m|f(2m+s)-f(2m+n)|\\
        &\leq A|b_m-c_m|+ C_{n,s}A\frac{c_m}{m+1}\\
        &\leq C_{n,s} A (m+1)^{s-2}
    \end{aligned}
    \]where the last inequality follows from \eqref{asintotici b e c} and \eqref{asintotico b-c}. Hence we obtain \[
    \left|\sum_{m=0}^{+\infty}\bigl[b_mf(2m+s)-c_mf(2m+n)\bigr]\right|\
 \le A C_{n,s}\sum _{m=0}^{+\infty}\frac{1}{(m+1)^{2-s}},
    \]and the last series is convergent being $s\in (0,1).$ From the estimate in \eqref{eq:cancel-small} arguing in a similar way we find \[
    |b_mf_h(2m+s)-c_mf_h(2m+n)|\leq C_{n,s,c}Ah(m+1)^{s-1}e^{-\tilde c(m+1)h}
    \]for some constant $0<\tilde c\leq c.$ From equation \eqref{eq:laplace-upper} in Lemma \ref{lem:laplace} we obtain \[
     \left|\sum_{m=0}^{+\infty}\bigl[b_mf_h(2m+s)-c_mf_h(2m+n)\bigr]\right|
 \le Ah\sum_{m=0}^{+\infty}C_{n,s,c}(m+1)^{s-1}e^{-\tilde c(m+1)h}\leq C_{n,s,c}Ah^{1-s}.
    \]
\end{proof}To estimate the derivatives of $R_i$ and $L_i$ we first recall an elementary estimate for alternating sums in terms of total variation.
We recall that a function $\psi:[x_1,+\infty)$ has finite total variation if \[\TV\psi:= \sup\sum_{j=1}^{l-1} \left|\psi(t_{j+1})-\psi(t_j)\right|<+\infty \] where the supremum is taken over all finite partitions \[ x_1\le t_1<\cdots<t_l<+\infty. \]It is well known that if $\psi \in C^1([x_1,+\infty))$ with $\psi'\in L^1([x_1,+\infty))$ then its total variation is equal to the $L^1$ norm of its derivative.
\begin{lemma}\label{lem:alt-bv}
Let $x_1<\cdots<x_K$ and let $\psi$ have bounded variation on
$[x_1,+\infty)$.  Then
\begin{equation}\label{eq:alt-bv}
 \left|\sum_{k=1}^K(-1)^{k+1}\psi(x_k)\right|
 \le |\psi(x_1)|+\TV_{[x_1,+\infty)}\psi.
\end{equation}
\end{lemma}
\begin{proof}
If $K$ is even we have \[
\left|\sum_{k=1}^{K}(-1)^{k+1}\psi(x_k)\right|\leq\sum_{k=1}^{K/2}|\psi(x_{2k-1})-\psi(x_{2k})|\leq \TV_{[x_1,+\infty)}\psi.
\]If $K$ is odd, observe instead
\[
 \left|\sum_{k=1}^K(-1)^{k+1}\psi(x_k)\right|
 \leq |\psi(x_1)|+\sum_{k=1}^{(K-1)/2}
 |\psi(x_{2k+1})-\psi(x_{2k})|\leq |\psi(x_1)|+\TV_{[x_1,+\infty)}\psi.
\]
\end{proof}
We can now use the previous lemma to obtain uniform estimates for the derivatives of the functions $R_i$, $L_i$ and $\Sigma_i$.
\begin{lemma}\label{lem:response-regularity}
There exists a constant $C>0$, such that for any finite equilibrium
\begin{equation}\label{eq:response-regularity}
 \lambda\bigl(|\partial_\lambda R_i(\lambda)|
 +|\partial_\lambda L_i(\lambda)|
 +|\partial_\lambda\Sigma_i(\lambda)|\bigr)\le C,
\end{equation}for any $\lambda>0.$
\end{lemma}
\begin{proof}
    For any $k=1,...,M-i$ we define $x_k=\sum_{q=i}^{i+k-1}p_q,$ in this way we have \[
    R_i(\lambda)=\sum_{k=1}^{M-i}(-1)^{k+1}e^{-\lambda x_k}.
    \]Differentiating with respect to $\lambda$ we find \[
    \partial_{\lambda}R_i(\lambda)=-\sum_{k=1}^{M-i}(-1)^{k+1}x_ke^{-\lambda x_k}.
    \]By explicit computation we can see that the maximum of the function $x\to x e^{-\lambda x}$ over $[0,+\infty)$ is $(e\lambda)^{-1}$, furthermore its total variation is $2/(e\lambda)$, hence from Lemma \ref{lem:alt-bv} we find \begin{equation}\label{prima derivata}
     |\partial_{\lambda}R_i(\lambda)|=\left|\sum_{k=1}^{M-i}(-1)^{k+1}x_ke^{-\lambda x_k}\right|\leq \frac{1}{e\lambda}+\frac{2}{e\lambda}=\frac{3}{e\lambda}.
    \end{equation}The estimate of $L_i$ is similar, indeed if we define $$y_k=\sum_{q=i-k}^{i-1}p_q,$$ we have \[
     L_i(\lambda)=\sum_{k=1}^{i-1}(-1)^{k+1}e^{-\lambda y_k}.
    \]Hence the same computations done in \eqref{prima derivata} show that \begin{equation}
     |\partial_{\lambda}L_i(\lambda)|\leq \frac{3}{e\lambda}\label{seconda derivata}.
    \end{equation}In order to estimate $\partial_\lambda \Sigma_i(\lambda)$ it is enough to combine the previous estimates with \eqref{eq:Sigma-div} and the first two inequalities in \eqref{eq:mode-bounds} to obtain \begin{equation}\label{terza derivata}
    |\partial_\lambda\Sigma_i(\lambda)|\leq|\partial_\lambda R_i(\lambda)||1-L_i(\lambda)|+|R_i(\lambda)||\partial_\lambda L_i(\lambda)|\leq \frac{6}{e\lambda}
    .\end{equation}
Combining \eqref{prima derivata}, \eqref{seconda derivata} and \eqref{terza derivata} we get \eqref{eq:response-regularity} with $C=12e^{-1}$.\end{proof}
For $1\le i\le M-1$ set
\begin{equation}\label{eq:rho-delta-zeta}
 \rho_i(\lambda):=\frac{e^{-\lambda p_i}}{1+e^{-\lambda p_i}},
 \qquad
 \delta_i(\lambda):=R_i(\lambda)-\rho_i(\lambda),
 \qquad
 \zeta_i(\lambda):=\rho_i(\lambda)\delta_i(\lambda).
\end{equation}From simple algebraic manipulation we obtain the following relations between the functions defined above.\begin{lemma}\label{lem:coboundary-mode}
For $1\le i\le M-2$ the following identity holds true 
\begin{equation}\label{eq:coboundary-mode}
 \delta_i(\lambda)
 =\zeta_i(\lambda)-\zeta_{i+1}(\lambda)
 +(\rho_{i+1}(\lambda)-\rho_i(\lambda))
  (\delta_{i+1}(\lambda)-\rho_i(\lambda)).
\end{equation}\end{lemma}\begin{proof}We observe that $e^{-\lambda p_i}=\rho_i(\lambda)(1-\rho_i(\lambda))^{-1},$ and $R_i(\lambda)=\rho_i(\lambda)+\delta_i(\lambda)$ hence, from \eqref{eq:recursions}, we find \[
\rho_i(\lambda)+\delta_i(\lambda)=e^{-\lambda p_i}(1-\rho_{i+1}(\lambda)-\delta_{i+1}(\lambda))=\frac{\rho_i(\lambda)}{1-\rho_i(\lambda)}(1-\rho_{i+1}(\lambda)-\delta_{i+1}(\lambda)),
\]from which we obtain \[
(1-\rho_i(\lambda))\delta_i(\lambda)=\rho_i(\lambda)(\rho_i(\lambda)-\rho_{i+1}(\lambda)-\delta_{i+1}(\lambda)).
\]From this identity and an explicit computation we conclude\[\begin{aligned}
\zeta_i(\lambda)-\zeta_{i+1}(\lambda)
 +(\rho_{i+1}(\lambda)-\rho_i(\lambda))
  (\delta_{i+1}(\lambda)-\rho_i(\lambda))&=\rho_i(\lambda)\delta_i(\lambda)+\rho_i(\rho_i-\rho_{i+1}-\delta_{i+1})\\&=\rho_i(\lambda)\delta_i(\lambda)+(1-\rho_i(\lambda))\delta_i(\lambda)\\&=\delta_i(\lambda).
\end{aligned}\]
\end{proof}
The above identities will allow us to approximate $\mathcal B_i-\Phi_{n,s}(e^{-p_i})$ with the following corrector terms.
For $0\le i\le M-1$ we introduce the following quantities
\begin{align}
 \mathcal C_i&:=\frac12\sum_m
 \bigl(b_m\Sigma_i({2m+s})-c_m\Sigma_i({2m+n})\bigr),
 \label{eq:Ccorr}\\
 \mathcal P_i&:=\sum_m
 \bigl(b_m\zeta_i({2m+s})-c_m\zeta_i({2m+n})\bigr),
 \qquad 1\le i\le M-1,
 \label{eq:Pcorr}\\
 \mathcal Q_i&:=-\mathcal C_i-\mathcal P_{i+1},
 \qquad 0\le i\le M-2.
 \label{eq:Qcorr}
\end{align}
Observe that $\mathcal C_0=0$ because $\Sigma_0(\lambda)=0$. As the following lemma states, the above quantities are well defined and bounded by a constant that depends only on $n$ and $s.$

\begin{lemma}\label{lem:corrector-bounded}
The series in \eqref{eq:Ccorr} and \eqref{eq:Pcorr} are absolutely convergent. Furthermore there exists a constant $C_{n,s}$ that is independent of $i$ and of the finite configuration such that
\begin{equation}\label{eq:Qbounded}
 |\mathcal C_i|+|\mathcal P_i|+|\mathcal Q_i|\le C_{n,s}.
\end{equation}
\end{lemma}
\begin{proof} By Lemmas~\ref{lem:mode} and~\ref{lem:response-regularity}, \[ |\Sigma_i(\lambda)| + \lambda|\partial_\lambda\Sigma_i(\lambda)| \le C \] for some universal constant $C$ that does not depend on the configuration. Hence \eqref{eq:cancel-bounded}, applied to \(f(\lambda)=\Sigma_i(\lambda)\), shows that the series defining \(\mathcal C_i\) converges and that $|\mathcal C_i|\le C_{n,s}.$ We next estimate \(\zeta_i\). Since \[ \rho_i(\lambda) = \frac{e^{-\lambda p_i}}{1+e^{-\lambda p_i}}, \] we have \[ 0\le\rho_i(\lambda)\le\frac12\quad\text{and}\quad \lambda|\partial_\lambda\rho_i(\lambda)| = \frac{\lambda p_i e^{-\lambda p_i}} {(1+e^{-\lambda p_i})^2} \le C. \] Moreover, \[ \delta_i(\lambda)=R_i(\lambda)-\rho_i(\lambda). \] Therefore, by Lemmas~\ref{lem:mode} and~\ref{lem:response-regularity}, \[ |\delta_i(\lambda)| + \lambda|\partial_\lambda\delta_i(\lambda)| \le C. \] Since $ \zeta_i(\lambda) = \rho_i(\lambda)\delta_i(\lambda)$, we easily infer \[ |\zeta_i(\lambda)| + \lambda|\partial_\lambda\zeta_i(\lambda)| \le C. \] A second application of \eqref{eq:cancel-bounded}, now with \(f(\lambda)=\zeta_i(\lambda)\), shows that the series defining \(\mathcal P_i\) converges and that \[ |\mathcal P_i|\le C_{n,s}. \] Finally, from \eqref{eq:Qcorr} we obtain \[ |\mathcal Q_i| \le |\mathcal C_i|+|\mathcal P_{i+1}| \le C_{n,s}, \]that implies \eqref{eq:Qbounded}. \end{proof}
We are finally ready to state and prove the main result of this subsection.
\begin{theorem}\label{thm:phase-corrector}
Fix $0<\vartheta<\Theta$.  
There is $C>0$ such that, whenever
$1\le i\le M-2$ and
\[
 \vartheta h\le p_i,p_{i+1}\le\Theta h,
 \qquad 
 0<h\leq1
\]
one has
\begin{equation}\label{eq:phase-corrector}
 \mathcal B_i-\Phi_{n,s}(e^{-p_i})
 =\mathcal Q_i-\mathcal Q_{i-1}+\mathcal E_i,
 \end{equation}for some error $
 |\mathcal E_i|\le Ch^{1-s}
$ and $C$ depends on $\vartheta$ and $\Theta.$
\end{theorem}
\begin{proof}
    From \eqref{eq:Sigma-div} and the definition of $\mathcal C _i$ we have \[
    \begin{aligned}
        \mathcal C_i-\mathcal C_{i-1}&=\frac{1}{2}\sum_{m=0}^{+\infty}b_m(\Sigma_i(2m+s)-\Sigma_{i-1}(2m+s))-c_m(\Sigma_i(2m+n)-\Sigma_{i-1}(2m+n))\\
        &=\frac{1}{2}\sum_{m=0}^{+\infty}b_m(R_i(2m+s)-L_i(2m+s))-c_m(R_i(2m+n)-L_i(2m+n)).
    \end{aligned}
    \]Hence, if we define \[
    \mathcal U_i=\sum_{m=0}^{+\infty}b_mR_i(2m+s)-c_mR_i(2m+n)
    \]we can write $\mathcal B_i=\mathcal U_i-(\mathcal C_i-\mathcal C_{i-1}).$
    Thanks to this definition and from \eqref{eq:rho-delta-zeta} we also obtain
    \[
    \mathcal U_i-\Phi_{n,s}(e^{-p_i})=\sum_{m=0}^{+\infty}b_m\delta_i(2m+s)-c_m\delta_i(2m+n),
    \]that thanks to \eqref{eq:coboundary-mode} can be rewritten as \[
     \begin{aligned} 
     \mathcal U_i-\Phi_{n,s}(e^{-p_i})&=\sum_{m=0}^{+\infty}b_m\zeta_i(2m+s)-c_m\zeta_i(2m+n)\\&\hspace{1em}-\sum_{m=0}^{+\infty}b_m\zeta_{i+1}(2m+s)-c_m\zeta_{i+1}(2m+n)+\mathcal E_i\\
     &=\mathcal P_i-\mathcal P_{i+1}+\mathcal E_{i}
   \end{aligned} \]where \[\begin{aligned}
   \mathcal E_i&=\sum_{m=0}^{+\infty}b_m(\rho_{i+1}(2m+s)-\rho_i(2m+s))(\delta_{i+1}(2m+s)-\rho_i(2m+s))\\
   &\hspace{1em}-\sum_{m=0}^{+\infty}c_m(\rho_{i+1}(2m+n)-\rho_i(2m+n))(\delta_{i+1}(2m+n)-\rho_i(2m+n)).
  \end{aligned} \]Thanks to the above equation and recalling that $\mathcal B_i=\mathcal U_i-(\mathcal C_i-\mathcal C_{i-1}),$ we infer \[
   \mathcal B_i-\Phi_{n,s}(e^{-p_i})
 =\mathcal Q_i-\mathcal Q_{i-1}+\mathcal E_i,
  \]where we also took into account \eqref{eq:Qcorr}. It remains only to show that $|\mathcal E_i|\leq C h^{1-s}$ for some constant $C.$ The assumption on $p_i$ and $p_{i+1}$ and a straightforward computation show that \[
   |\rho_{i+1}(\lambda)-\rho_i(\lambda)|
 \le C\lambda h e^{-c\lambda h}\quad\text{and}\quad
 |\partial_\lambda(\rho_{i+1}(\lambda)-\rho_i(\lambda))|
 \le Ch(1+\lambda h)e^{-c\lambda h},
  \]moreover, in the proof of Lemma \ref{lem:corrector-bounded} we showed that the functions $\rho_i(\lambda)$ and $\delta_i(\lambda)$ are bounded and their derivatives fulfill $\lambda|\partial_\lambda\rho_i(\lambda)|\leq C$ and $\lambda|\partial_\lambda\delta_i(\lambda)|\leq C.$ From this observation we infer \[
  |(\rho_{i+1}(\lambda)-\rho_i(\lambda))(\delta_{i+1}(\lambda)-\rho_i(\lambda))|\leq C\lambda h e^{-c\lambda h}
  \]and \[
   |\partial_{\lambda}[(\rho_{i+1}(\lambda)-\rho_i(\lambda))(\delta_{i+1}(\lambda)-\rho_i(\lambda))]|\leq Ch(1+\lambda h)  e^{-c\lambda h},
  \]that together with \eqref{eq:cancel-small} implies that $\mathcal |E_i|\leq C h^{1-s}.$
\end{proof}
\subsection{Polynomial decay of the first radius}\label{sec:packing}
In this subsection we prove the last part of Theorem \ref{thm:main finite}, that is in any finite equilibrium the first radius goes to zero polynomial in the number of shells. \eqref{eq:packing-main} holds. We start by estimating the function $\mathcal M_i$ when three adjacent gaps have the same scale.
\begin{lemma}\label{lem:stress-comparison}The following facts hold true:\begin{enumerate}[label=(\roman*)]
\item There exists a constant $C_{n,s}$ such that for every $1\le i\le M-1$,
\begin{equation}\label{eq:stress-upper-general}
 \mathcal M_i\le C_{n,s}p_i^{-s}.
\end{equation}
\item For every $0<\theta<\Theta$ there are $c,C$ that depend only on $n$, $s,$ $\theta$, and $\Theta$ such that for every
$2\le i\le M-2$,
\begin{equation}\label{eq:stress-comparison}
 ch^{-s}\le\mathcal M_i\le Ch^{-s}
\end{equation}
provided that $\theta h\le p_{i-1},p_i,p_{i+1}\le\Theta h$ and $0<h<\frac14$.
\end{enumerate}
\end{lemma}
\begin{proof}
    Let us start with the proof of $(i).$ From the third estimate in \eqref{eq:mode-bounds} we have $\Sigma_i(\lambda)\leq e^{-\lambda p_i}$, hence from \eqref{eq:Mstress} and \eqref{asintotici b e c} we find \[
    \mathcal M_i\leq \frac{1}{2}\sum_{m=0}^{+\infty}b_m e^{-(2m+s)p_i}+c_me^{-(2m+n)p_i}\leq C_{n,s}\sum_{m=0}^{+\infty}(m+1)^{s-1}e^{-c_{n,s}(m+1)p_i}.
    \]Now we distinguish two cases. If $p_i\leq1$ we apply \eqref{eq:laplace-upper} in Lemma \ref{lem:laplace} to obtain \[
    \sum_{m=0}^{+\infty}(m+1)^{s-1}e^{-c_{n,s}(m+1)p_i}\leq C_{n,s}p_i^{-s}.
    \]If $p_i>1$ we simply observe that \[
    \sum_{m=0}^{+\infty}(m+1)^{s-1}e^{-c_{n,s}(m+1)p_i}\leq e^{-c_{n,s}p_i}\sum_{m=0}^{+\infty}(m+1)^{s-1}e^{-c_{n,s}m}\leq C_{n,s}e^{-c_{n,s}p_i}\leq C_{n,s}p_i^{-s},
    \]hence in both cases we have \eqref{eq:stress-upper-general}.
    
    We now deal with $(ii)$, the upper bound simply follows by \eqref{eq:stress-upper-general} and the assumption on $p_i$. We now deal with the lower bound, by \eqref{eq:recursions}, \eqref{eq:Sigma-div} and \eqref{eq:mode-bounds} we infer \[
    \Sigma_i(2m+s)\geq e^{-(2m+s)p_i}(1-e^{-(2m+s)p_{i+1}})(1-e^{-(2m+s)p_{i-1}}).
    \]If we choose $m$, depending on $h$ and $s$ such that $1\leq h(2m+s)\leq 2$ we have $(2m+s)p_i\leq 2\Theta$, $(2m+s)p_{i+1}\geq\theta$ and $(2m+s)p_{i-1}\geq\theta$, hence we find \[
    \Sigma_i(2m+s)\geq e^{-2\Theta}(1-e^{-\theta})^2=C_{\theta,\Theta}.
    \]So, recalling also \eqref{asintotici b e c}, we obtain \begin{equation}\label{stima basso}
    \mathcal M_i\geq\frac12\sum_{m=0}^{+\infty}b_m\Sigma_i(2m+s)\geq {C_{\theta,\Theta,n,s}} \sum_{1\leq (2m+s)h\leq 2} (m+1)^{s-1}.
    \end{equation}It is easy to see that, for any $h<\frac14$, any integer that fulfills $$\frac{3}{4h}\leq m+1\leq \frac{1}{h}$$ also fulfills $1\leq (2m+s)h\leq2.$ Hence from \eqref{stima basso} we conclude \[
    \mathcal M_i\geq C_{n,s,\theta,\Theta}\sum_{\frac{3}{4h}\leq m+1\leq \frac{1}{h}}(m+1)^{s-1}\geq C_{n,s,\theta,\Theta} h^{-s},
    \]where the last estimate easily follows by integral comparison.
\end{proof}
The following lemma, whose proof immediately follows from \eqref{eq:adjacent-main} implies that two consecutive gaps cannot skip the band $[\lambda h,h].$
\begin{lemma}\label{lem:no-skip}
Let $c_0,C_0,\eta_0$ be the constants in Theorem~\ref{thm:adjacent}, and fix
\[
 0<\lambda<
 \min\{c_0,C_0^{-1}\}.
\]Then, for every 
$h>0$ and every index $1\leq j\leq M-3$  with $\min\{p_{j},p_{j+1}\}\leq \eta_0$ the following facts hold: \begin{enumerate}[label=(\roman*)]
    \item If $p_{j}> h$ then $p_{j+1}>\lambda h$.
    \item If $p_{j+1}>h$ then $p_j>\lambda h.$
\end{enumerate}
\end{lemma}
\begin{proof}Both estimates follow from \eqref{eq:adjacent-main} in Theorem \ref{thm:main finite}. Indeed if $p_j>h$ then, from $\eqref{eq:adjacent-main}$, $p_{j+1}\geq c_0 p_j>h\lambda.$ In the same way, if $p_{j+1}>h$ then, again from \eqref{eq:adjacent-main}, we infer $p_j\geq C_0^{-1}p_{j+1}> h\lambda.$ 
\end{proof}
The following definition describes a transition from the scale $h$ to the smaller scale $\lambda h.$
\begin{definition}\label{def:crossing}
Fix $0<\lambda<1$.  For indices $2\le a<b\le M-1$, the pair
$(a,b)$ is a downward crossing of the band
$[\lambda h,h]$ if $p_{a-1}>h,$ $p_b\le\lambda h$ and
\begin{equation}\label{eq:crossing}
 \lambda h<p_j\le h\quad\text{for all $a\leq j<b$,}
\end{equation}
\end{definition}
We now show that the sum of the logarithmic gaps along every downward crossing is bounded from below by a fixed positive constant.
\begin{theorem}\label{thm:transition}
    There exists $\lambda_0=\lambda_0(n,s)\in(0,1)$ and a constant $\tilde\eta=\tilde\eta(n,s,\lambda_0)>0$ such that, for every \(0<\lambda<\lambda_0\), there exists $\delta_*=\delta_*(n,s,\lambda)>0$ with the following property. Every downward crossing \((a,b)\) of the band $[\lambda h,h],$ as in Definition \ref{def:crossing}, with $0<h\le\tilde\eta $ {and} $ b\le M-3,$ satisfies\begin{equation}\label{eq:transition-cost}
\log\bigg(\frac{r_b}{r_a}\bigg)= \sum_{j=a}^{b-1}p_j\ge\delta_*.
\end{equation}\end{theorem}
\begin{proof}
    Let $0<\lambda<\lambda_0$, where $\lambda_0$ will be chosen later, but small enough that Lemma \ref{lem:no-skip} holds true. Let  $(a,b)$ be a downward crossing of the band $[\lambda h,h]$ with $\lambda<\lambda_0$ . If we choose $\tilde\eta\leq \min\{\eta_0,\frac14\}$ where $\eta_0$ is the one in \eqref{eq:adjacent-main} of Theorem \ref{thm:main finite} we find \[
    c_0 h<c_0p_{a-1}\leq p_a\leq h.
    \]Hence from \eqref{eq:stress-upper-general} we find, up to a possibly different constant, \begin{equation}\label{sistemo il primo}\mathcal M_a\leq C_{n,s}p_a^{-s}\leq C_{n,s}h^{-s}.\end{equation}Now we deal with $p_b$. By definition we have $p_{b-1}>\lambda h$ and $p_b\leq \lambda h< \tilde\eta$, hence by \eqref{eq:adjacent-main} we infer $c_0p_{b-1}\leq p_b\leq C_0p_{b-1},$ that implies \begin{equation}
    \lambda h< p_{b-1}\leq c_0^{-1}\lambda h \quad\text{and}\quad c_0\lambda h\leq p_b\leq\lambda h.
    \end{equation}From this estimate and applying again \eqref{eq:adjacent-main}\begin{equation}\label{caso speciale}
    c_0^2\lambda h\leq c_0p_b\leq p_{b+1}\leq C_0p_b\leq C_0\lambda h.
    \end{equation}Hence we can apply \eqref{eq:stress-comparison} to the jumps $p_{b-1}$, $p_b$ and $p_{b+1}$ with scale $\lambda h$ to infer that \begin{equation}\label{sistemo l'ultimo}\mathcal M_b\geq c (\lambda h)^{-s}.\end{equation}
    Combining \eqref{sistemo il primo} with \eqref{sistemo l'ultimo} and choosing $\lambda_0 $  small enough we find \[
    \mathcal M_b-\mathcal M_a\geq c (\lambda h)^{-s}-C_{n,s}h^{-s}\geq  c (\lambda_0 h)^{-s}-C_{n,s}h^{-s}\geq ch^{-s}
    \]for a constant that depends on $\lambda_0$. Now, from Lemma \ref{lem:corrector-bounded} we know that there exists a constant $C_{n,s}$ such that $|\mathcal Q_a|\leq C_{n,s}$ and $|\mathcal Q_b|\leq C_{n,s},$ from this we find, up to choosing $\tilde\eta$ small enough \begin{equation}\label{primo e ultimo}
    (\mathcal M_b+\mathcal Q_b)-(\mathcal M_a+\mathcal Q_a)\geq Ch^{-s}-C_{n,s}\geq c h^{-s}.
    \end{equation}
   Now fix $a<i\leq b,$ since the couple $(a,b)$ is a downward crossing (and by \eqref{caso speciale} for $i=b$) we infer the existence of constants $\theta$ and $\Theta$ (depending on $n$, $s$ and $\lambda$) such that \[
   \theta h\leq p_i,p_{i+1}\leq \Theta h\quad \text{for }a<i\leq b.
   \]Hence we can apply Theorem \ref{thm:phase-corrector} to infer \[
   \mathcal B_i-\Phi_{n,s}(e^{-p_i})=\mathcal Q_i-\mathcal Q_{i-1}+\mathcal E_i,
   \]that combined with \eqref{eq:stress-balance} implies, after some algebraic manipulations \begin{equation}\label{parti in mezzo}
   (\mathcal M_i+\mathcal Q_i)-(\mathcal M_{i-1}+\mathcal Q_{i-1})=\Phi_{n,s}(1^-)-\Phi_{n,s}(e^{-p_i})-\mathcal E_i.
   \end{equation}From \eqref{eq:defect-two-sided}, the fact that $p_i\leq h$, and by choosing $\tilde\eta\leq p_0$ we infer\footnote{We are also using that the function $p\to\Phi_{n,s}(1^-)-\Phi_{n,s}(e^{-p})$ is increasing, this follows from Lemma \ref{lem:Phi}.} \[
   \Phi_{n,s}(1^-)-\Phi_{n,s}(e^{-p_i})\leq C h^{1-s}.
   \]Combining this estimate with \eqref{parti in mezzo} and recalling that, from Theorem \ref{thm:phase-corrector}, $|\mathcal E_i|\leq C h^{1-s}$ (for some constant $C$ that depends also on $\lambda$) we find \begin{equation}
       \label{sistemo quelli in mezzo}|(\mathcal M_i+\mathcal Q_i)-(\mathcal M_{i-1}+\mathcal Q_{i-1})|\leq C h^{1-s}.
   \end{equation}
   From \eqref{primo e ultimo} and \eqref{sistemo quelli in mezzo} we find \[
   c h^{-s}\leq (\mathcal M_b+\mathcal Q_b)-(\mathcal M_a+\mathcal Q_a)\leq \sum_{i=a+1}^b|(\mathcal M_i+\mathcal Q_i)-(\mathcal M_{i-1}+\mathcal Q_{i-1})|
   \leq C(b-a)h^{1-s},\]that implies $(b-a)\geq \tilde c h^{-1}$ for a suitable constant $\tilde c$ that depends on $n$, $s$ and $\lambda.$
Since for any $a\leq j<b$ it holds $p_j\geq \lambda h$ we finally conclude \[
\sum_{j=a}^{b-1}p_j\geq \lambda h(b-a)\geq \lambda \tilde c:=\delta_*
,\]where $\delta_*$ depends also on $\lambda.$\end{proof}
We now estimate the number of disjoint downward crossings between two fixed scales.
\begin{lemma}\label{lem:crossing-extraction}
Let $q_{j_0-1},q_{j_0},\ldots,q_L>0$, with
$q_{j_0-1}>\eta$ and $q_{j_0}\in[c\eta,\eta]$ with $0<c<1$.  Suppose
\begin{equation}\label{assunzione extraction}\min_{j_0\le j\le L}q_j\le\mu<\eta,\end{equation}and assume the existence of constants $0<c_0<C_0$ such that \begin{equation}\label{adjacent lemma}
\min\{q_j,q_{j+1}\}\leq \eta\quad\implies\quad c_0 q_j\leq q_{j+1}\leq C_0 q_j,\quad\text{if }\, j_0-1\leq j<L
.\end{equation}Fix also $ 0<\lambda<\min\{c_0,C_0^{-1}\}. $ Then one can find at least \begin{equation}\label{numero salti}
\frac{\log(\eta/\mu)}{|\log(\lambda)|}-1
\end{equation}downward crossings $(a_k,b_k)$ with $j_0\leq a_k<b_k\leq L$ and associated bands $[\lambda^{k+1}\eta,\lambda^k\eta ]$. Furthermore the intervals $[a_k,b_{k}-1]$ are pairwise disjoint.\end{lemma}
\begin{proof}
    For any $h_k:=\lambda^k\eta.$  If $k$ is such that $h_{k+1}\geq\mu $, let $\tau_k$ be the first index $j$ such that $q_j\leq h_k,$ namely \[
    \tau_k:=\min\big\{{j\in\{j_0,...,L\}}\::\: q_j\leq h_k\big\}.
    \]We now show that the sequence $\tau_k$ is increasing. Indeed since $q_{\tau_k}\leq h_k<\eta$ we can apply \eqref{adjacent lemma} to obtain \[
    q_{\tau_k}\geq c_0 q_{\tau_k-1}>c_0 h_k>\lambda h_k=h_{k+1}
    ,\]that implies $\tau_{k}<\tau_{k+1}.$ Let now $b_k=\tau_{k+1}$ and let $a_k-1$ be the largest index that is smaller than $b_k$ such that $q_{a_k-1}>h_k.$ 
    We notice that $a_k<b_k$, indeed if we would have $a_k=b_k$ that implies $q_{b_k-1}>h_k$ and $q_{b_k}\leq \lambda h_k$, with a similar argument as the one of Lemma \ref{lem:no-skip} we find a contradiction.
    For this choice of $a_k$ we have \[
    q_j\leq h_k\quad\text{for all $a_k\leq j<b_k$}
    \]and, recalling $b_k=\tau_{k+1}$ is the first index $j$ such that $q_j\leq h_{k+1}$ we infer \[
q_j>h_{k+1}\quad\text{for all $a_k\leq j<b_k$.}
    \]So we have \[
    h_{k+1}<q_j\leq h_k\quad\text{for all $a_k\leq j<b_k$},
    \]that combined with the fact that $q_{a_k-1}>h_k$ and $q_{b_k}\leq h_{k+1}$ shows that $(a_k,b_k)$ is a downward crossing of the band $[h_{k+1},h_k]$. By construction we have $a_{k+1}\geq \tau_{k+1}=b_k$ which implies that the intervals $[a_k,b_k-1]$ are pairwise disjoint. It remains only to prove that the number of crossings is at least \eqref{numero salti}. In order to do this it is enough to observe that the construction of the downward crossing $[h_k,h_{k+1}]$ is possible if $h_{k+1}=\lambda^{k+1}\eta\geq\mu.$ Solving for $k$ gives \eqref{numero salti}.
\end{proof}
We can now use the cost of the disjoint crossings to obtain the polynomial decay of the inner radius.
\begin{theorem}\label{thm:packing}
There are $C_{n,s},\sigma_{n,s}>0$ such that every finite radial equilibrium satisfies
\eqref{eq:packing-main}.
\end{theorem}
\begin{proof}
    We fix:\begin{itemize} 
    \item $\lambda<\min\{\lambda_0,c_0,C_0^{-1}\}$ where $\lambda_0$ is the one of Theorem \ref{thm:transition} and $c_0$, $C_0$ are the one of \eqref{eq:adjacent-main};
    \item $\eta<\min\{\tilde\eta,\,\eta_0,\log(1+c_{n,s})\}$ where $\tilde\eta$ is the one of Theorem \ref{thm:transition}, $\eta_0$ is the one in \eqref{eq:adjacent-main} and $c_{n,s}$ is the one of Theorem \ref{thm:central}.
    \end{itemize}
After enlarging the final constant, we may assume that $M\geq 12$. If $-\log(r_1)\geq\frac{\eta M}{4}$ we have $-\log(r_1)\geq\frac\eta4\log(M)$ and hence \[
    r_1=e^{\log(r_1)}\leq M^{-\frac \eta 4}.
    \]We now assume $-\log(r_1)<\frac{\eta M}{4}.$ We claim that there exists at least one index $1\leq j\leq M-3$ such that $p_j\leq \eta.$ Indeed, if this is not the case, we would have \[
    -\log(r_1)=\sum_{j=1}^{M-1}p_j\geq \sum_{j=1}^{M-3}p_j\geq (M-3)\eta\geq \frac{\eta M}4
    \]where the last inequality follows from the fact that $M\geq 12.$ Let $j_0$ be the first index such that $p_{j_0}\leq \eta,$ since $p_1>\eta$ we have $j_0\geq 2$ and $p_j>\eta$ for any $1\leq j<j_0$. So we have \[
    (j_0-1)\eta<\sum_{j=1}^{j_0-1}p_j\leq -\log(r_1)<\frac{M\eta}{4}
    \]that implies $j_0-1<\frac{M}{4}.$ Hence we can apply Theorem \ref{thm:adjacent} to infer, recalling also that $p_{j_0}\leq \eta$ \begin{equation}\label{prima ipotesi}
    c_0\eta<c_0 p_{j_0-1}\leq p_{j_0}\leq\eta.
    \end{equation}Let us now consider the indices $j_0,j_0+1,...,M-3,$ since $j_0-1<\frac M4$ and $M\geq 12$ we can see that the set $\{p_{j_0},...,p_{M-3}\}$ contains at least $\frac M2$ elements. Since their sum is at most $-\log r_1$ there is at least one of them, let us say $p_j$ such that $$p_j\leq \frac{-2\log(r_1)}{M}<\frac\eta2<\eta.$$Hence \eqref{assunzione extraction} of Lemma \ref{lem:crossing-extraction} holds true with the choice $\mu=\frac{-2\log(r_1)}{M}.$ Therefore by Lemma \ref{lem:crossing-extraction} we find at least \[
    \widetilde N\geq\frac{1}{|\log(\lambda)|}\log\bigg(\frac{\eta M}{-2\log(r_1)}\bigg)-1
    \]
    downward crossings $(a_k,b_k)$ with disjoint corresponding intervals. Hence from Theorem \ref{thm:transition} we find \[
    -\log(r_1)=\sum_{j=1}^{M-1}p_j\geq \sum_{k=1}^{\widetilde N}\sum_{j=a_k}^{b_{k}-1}p_j\geq \sum_{k=1}^{\widetilde N}\delta_*=\widetilde N\delta_*.
    \]That implies \begin{equation}\label{quasi finita}
    -\log(r_1)\geq a\log\bigg(\frac{cM}{-\log(r_1)}\bigg)-C
    \end{equation}for suitable constant $a$, $c$ and $C.$ We observe that from \eqref{quasi finita} the required estimates easily follows, indeed if $-\log r_1\geq\sqrt M$ \eqref{eq:packing-main} is trivial. Otherwise we observe that \[
    \log\bigg(\frac{cM}{-\log(r_1)}\bigg)\geq \frac{1}{2}\log(M)+\log(c)
    \]that combined with \eqref{quasi finita} implies \[
    -\log(r_1)\geq\frac{a}{2}\log (M)+a\log(c)-C
    \]that is equivalent to \eqref{eq:packing-main}.
\end{proof}

\section{Lower bound and limit set}\label{sec:lower-bound-limit}

In this last section, we finally prove a uniform lower bound for the gaps, which will allow us to eventually take the limit as the number of shells goes to $+\infty$.
As in the end of Section \ref{sec:existence-finite-equilibria}, for each even $M$ choose one maximizer from Theorem~\ref{thm:finite-existence}, denote its radii as $r_i^{(M)}$ and then define the gaps
\begin{equation*}
 p_i^{(M)}=\log\frac{r_{i+1}^{(M)}}{r_i^{(M)}}.
\end{equation*}

\begin{lemma}\label{lem:fixed-index-compactness}
There are numbers $d_i>0$, depending only on $n$, $s$, and $i$, such that for every fixed $i$ and every sufficiently large even $M$,
\begin{equation*}
 d_i\le p_i^{(M)}\le\overline p_{n,s}.
\end{equation*}
In particular, by compactness we can deduce that, up to choosing a subsequence, which we won't rename for the sake of readability, there exist values $p_i>0$ such that
\begin{equation}\label{eq:gap-diagonal-limit}
 p_i^{(M)}\longrightarrow p_i \qquad\text{for every fixed }i.
\end{equation}
\end{lemma}

\begin{proof}
The upper bound has already been proved in Lemma~\ref{lem:gap-upper-all}.
To prove the existence of the $d_i$'s we proceed by induction.
Theorem~\ref{thm:central} gives $p_1^{(M)}\ge d_1$ for a constant $d_1=d_1(n,s)>0$.  
Let now $c_0,\eta_0$ be the lower comparison constant and threshold in Theorem~\ref{thm:adjacent}, and define recursively
\[
 d_{i+1}:=\min\{\eta_0,c_0d_i\}.
\]
Assume $p_i^{(M)}\ge d_i$.  If $p_{i+1}^{(M)}\ge\eta_0$, then $p_{i+1}^{(M)}\ge d_{i+1}$.  
Otherwise, the hypothesis of Theorem~\ref{thm:adjacent} holds and 
\[p_{i+1}^{(M)}\ge c_0p_i^{(M)}\ge c_0d_i\ge d_{i+1}.\]
For any fixed $i$, all of the required adjacent comparisons are available once $M\ge i+3$.  
This proves the lower bound and a diagonal extraction gives \eqref{eq:gap-diagonal-limit}.
\end{proof}

Now we are ready to rescale the finite equilibria $E_M$ and consider the limit as $M\to +\infty$.
Let us set some notation. 
Define
\begin{equation}\label{eq:rescaled-radii}
 \rho_i^{(M)}:=\frac{r_i^{(M)}}{r_1^{(M)}},
 \qquad
 \widetilde E_M:=\frac1{r_1^{(M)}}E(\mathbf r^{(M)}).
\end{equation} 
A quick change of variables shows that 
\begin{equation}\label{eq:H-scaling}
 H_s(\lambda x,\lambda E)=\lambda^{-s}H_s(x,E).
\end{equation}
for all $\lambda>0$ and $E\subset \R^n$, hence $\widetilde E_M$ is an equilibrium on its first $M-1$ interfaces.  
In particular, Theorem~\ref{thm:packing} gives
\begin{equation}\label{eq:outer-radius-diverges}
 \rho_M^{(M)}=\frac{1}{r_1^{(M)}}\ge C_{n,s}^{-1}M^{\sigma_{n,s}}\longrightarrow+\infty.
\end{equation}

Now, define
\begin{equation}\label{eq:limit-radii}
 \rho_1:=1,
 \qquad
 \rho_i:=\exp\bigg(\sum_{j=1}^{i-1}p_j\bigg),
 \quad i\ge2.
\end{equation}
Then, for every fixed $i$, since for all $j< i$ we have that $p_j^{(M)}\to p_j$ as $M\to+\infty$, we deduce
\begin{equation}\label{eq:radii-diagonal-limit}
 \rho_i^{(M)}\underset{M}{\longrightarrow}\rho_i,
 \qquad
 1=\rho_1<\rho_2<\cdots.
\end{equation}

The aim of the rest of the section is to prove that the limit set is a union of annuli with the $\rho_i$'s as radii and that the expression for the curvature is stable for $M\to +\infty$, thus deducing that the limit set is an $s$-minimal surface.
The first step in this direction is proving that $\rho_i$'s diverge as $i\to+\infty$.
Indeed, the fixed-index lower bounds in Lemma~\ref{lem:fixed-index-compactness} may deteriorate with $i$ and do not by themselves imply $\rho_i\to+\infty$. 
\begin{proposition}\label{prop:limit-radii-diverge}
The sequence in \eqref{eq:limit-radii} satisfies
\begin{equation*}
 \sum_{i=1}^{+\infty} p_i=+\infty,
 \qquad\text{and hence}\qquad
 \rho_i\longrightarrow +\infty.
\end{equation*}
\end{proposition}

\begin{proof}
Suppose by contradiction that $\sum_i p_i<+\infty$.  
Then, $p_i\to0$.
Pass the adjacent comparison inequalities to the limit on every fixed pair: whenever $\min\{p_i,p_{i+1}\}<\eta_0$, the strict inequality persists for
large $M$, and
\begin{equation}\label{eq:limit-adjacent-comparison}
 c_0p_i\le p_{i+1}\le C_0p_i.
\end{equation}

Choose the constants $\tilde \eta,\lambda,\delta_*$ given by Theorem~\ref{thm:transition}.  Choose a band top $h_0<\min\{\tilde\eta,p_1,\eta_0\}$ outside the countable set
\[
 \{\lambda^{-k}p_j:j\ge1,\ k\ge0\};
\]
then, no $p_j$ equals a boundary of any band $[\lambda^{k+1}h_0,\lambda^kh_0]$.

Let $j_0$ be the first index with $p_{j_0}<h_0$. 
The relation \eqref{eq:limit-adjacent-comparison} gives $p_{j_0}\ge c_0h_0$. 
Now, fix $K\ge1$.  Since $p_i\to0$, we can choose a finite $L=L(K)$ for which
\[
 \min_{j_0\le j\le L}p_j
\]
is small enough that Lemma~\ref{lem:crossing-extraction} yields at least $K$ pairwise disjoint downward crossings $(a_k,b_k)$, $1\le k\le K$, of bands $[\lambda h,h]$ with $h$ in the geometric family $\{\lambda^qh_0:q\ge0\}$. 
All the inequalities defining these crossings are strict thanks to the choice of $h_0$.

These $K$ crossings involve only finitely many indices.  For large $M$, convergence \eqref{eq:gap-diagonal-limit} preserves all their defining inequalities, and $M>\max_{k\le K}b_k+3$.  
Thus, Theorem~\ref{thm:transition} applies to the corresponding crossings in the finite equilibrium and gives
\[
 \sum_{j=a_k}^{b_k-1}p_j^{(M)}\ge\delta_*,
 \qquad k=1,\ldots,K.
\]
Passing to the limit in each finite sum then yields
\[
 \sum_{j=a_k}^{b_k-1}p_j\ge\delta_*.
\]
The crossing interiors are disjoint, so
\[
 \sum_{j=1}^{+\infty} p_j\ge K\delta_*.
\]
Since $K$ is arbitrary, this contradicts the assumed finite total mass.
Therefore, the conclusion follows.
\end{proof}

We're finally ready to describe the limit set and prove that it satisfies the equilibrium equation at every point of its boundary.

Set
\begin{equation}\label{eq:limit-set}
 E_\infty:=\bigcup_{k=1}^{+\infty} \bigl(B_{\rho_{2k}}\setminus B_{\rho_{2k-1}}\bigr).
\end{equation}

\begin{remark}
By Proposition~\ref{prop:limit-radii-diverge}, every bounded region meets only
finitely many of the spheres $\partial B_{\rho_i}$.  To verify the local
perimeter assertion explicitly, fix $L>0$ and choose $K$ with $\rho_K>2L$.
Inside $B_{2L}$ the set is a finite union of smooth annuli, so all interactions
with both variables in $B_{2L}$ are finite.  Interactions with one variable in
$B_L$ and the other outside $B_{2L}$ are bounded by
\[
 C_{n,s}|B_L|\int_L^{+\infty} \frac{1}{r^{1+s}}\,\dd r<+\infty.
\]
Hence, $E_\infty$ has locally finite fractional perimeter and a smooth boundary
near every boundary point.    
\end{remark}

\begin{theorem}\label{thm:unbounded-stationary}
The set $E_\infty$ in \eqref{eq:limit-set} satisfies
\begin{equation}\label{eq:limit-curvature-zero}
 H_s(x,E_\infty)=0
 \qquad\text{for every }x\in\partial E_\infty.
\end{equation}
\end{theorem}

\begin{proof}
Fix $i\ge1$ and define
\[
 x_M:=\rho_i^{(M)}e_1,
 \qquad x:=\rho_ie_1.
\]
Then $x_M\to x$.  Write
\[
 \sigma_M:=\one_{\widetilde E_{M}^c}
 -\one_{\widetilde E_{M}},
 \qquad
 \sigma_\infty:=\one_{E_\infty^c}-\one_{E_\infty},
\]
and introduce the reference phases
\begin{align*}
 \tau_M(y)&:=(-1)^i  \bigl(\one_{B_{\rho_i^{(M)}}^c}(y) -\one_{B_{\rho_i^{(M)}}}(y)\bigr),\\
 \tau_\infty(y)&:=(-1)^i \bigl(\one_{B_{\rho_i}^c}(y)-\one_{B_{\rho_i}}(y)\bigr).
\end{align*}
The two adjacent limiting gaps are positive.  Since $\rho_1^{(M)}=1$ for all $M$, for $i=1$ only the outer gap is needed.
Hence, there is $d_i>0$ such that, for large $M$, no interface of $\widetilde E_{M}$ other than $\partial B_{\rho_i^{(M)}}$ meets the radial annulus
\[
 \left\{y:\left||y|-\rho_i^{(M)}\right|<d_i\right\}.
\]
On this annulus $\sigma_M=\tau_M$.
The principal value can therefore be written as the absolutely convergent expression
\begin{equation}\label{eq:reference-phase-curvature}
 H_s(x_M,\widetilde E_{M}) =(-1)^i\frac{\kappaN}{\left(\rho_i^{(M)}\right)^s}  +\int_{\R^n}  \frac{\sigma_M (y)-\tau_M (y)}{|x_M-y|^{n+s}}\,\dd y.
\end{equation}

Now, fix $L>2\rho_i+2$.
On $B_L$, the local finiteness of the limiting radii and \eqref{eq:radii-diagonal-limit} imply $\sigma_M\to\sigma_\infty$ almost everywhere. 
Likewise, $\tau_M\to\tau_\infty$ almost everywhere.
The numerator in \eqref{eq:reference-phase-curvature} vanishes in a fixed neighborhood of $x_M$, while outside that neighborhood the kernel is uniformly bounded on $B_L$.  
The dominated convergence theorem therefore gives convergence of the integral on $B_L$.

For the tail, $|x_M|\le\rho_i+1$ for $M$ large enough, and
\begin{align}
 \int_{\R^n\setminus B_L} \frac{|\sigma_M (y)-\tau_M (y)|} {|x_M-y|^{n+s}}\,\dd y  \le C_{n,s}\int_L^{+\infty} \frac{1}{r^{1+s}}\,\dd r \le C_{n,s}L^{-s}.
\end{align}
The same estimate holds for the limiting phases. 
First, let $M\to+\infty$ at fixed $L$ in \eqref{eq:reference-phase-curvature}, then let $L\to+\infty$.
Since by scaling we do not alter the equilibrium equations, we have
\[
 H_s(x_M,\widetilde E_{M})=0.
\]
We conclude that $H_s(x,E_\infty)=0$.
By rotational invariance, the same identity holds at every point of $\partial B_{\rho_i}$, and since $i$ was arbitrary we conclude that $H_s(x,E_\infty)=0$ for all $x\in \partial E_{\infty}$ as desired.
\end{proof}

Theorem~\ref{thm:finite-existence} and
Theorem~\ref{thm:unbounded-stationary}, together with the quantitative results
of the preceding sections, complete the proof of Theorem~\ref{thm:main}.

\subsection*{Acknowledgments}
The authors are members of the Gruppo Nazionale per l’Analisi Matematica, la Probabilità e le loro Applicazioni (GNAMPA) of the Istituto Nazionale di Alta Matematica (INdAM). M.T. was partially funded by the INdAM–GNAMPA project "Problemi di ottimizzazione di forma in contesti anisotropi e non-locali", CUP: E53C25002010001.
C.G. was supported by the INdAM - GNAMPA project "Ottimizzazione geometrica e frontiere libere", CUP E53C25002010001.
The authors are grateful to Matteo Cozzi for stimulating discussions on this problem.
\bibliographystyle{abbrv} 
\bibliography{riferimenti}
\end{document}